\documentclass[fontsize=12pt,a4paper,headings=normal,
twoside=false,leqno,parskip=half-,abstract=true]{scrartcl}
\usepackage[english]{babel}
\usepackage[utf8]{inputenc}
\usepackage{hyperref}
\hypersetup{
 pdftitle={meanders with three noses},
 pdfauthor={Bernold Fiedler},
 colorlinks=true,
 linkcolor=blue,
 citecolor=blue,
 filecolor=blue,
 urlcolor=blue}

 \usepackage{afterpage}

\usepackage{graphicx}
\usepackage{wrapfig} 
\usepackage[format=plain,labelfont=bf,font=small]{caption}
\usepackage{xcolor}
\usepackage[arrow, matrix, curve]{xy}
\usepackage{float}

\usepackage{caption}
\usepackage{tabulary}
\usepackage{array}
\newcolumntype{N}[1]{>{\centering\arraybackslash}m{#1}}

\usepackage{amsmath,amsthm}
\swapnumbers 
\usepackage{amssymb, eurosym} 

\makeatletter
\newcommand{\tpitchfork}{%
  \vbox{
    \baselineskip\z@skip
    \lineskip-.52ex
    \lineskiplimit\maxdimen
    \m@th
    \ialign{##\crcr\hidewidth\smash{$-$}\hidewidth\crcr$\pitchfork$\crcr}
  }%
}
\makeatother
\usepackage{latexsym}
\usepackage{enumerate}

\usepackage[notref,notcite,color, final 
]{showkeys}

\definecolor{refkey}{rgb}{1,0,0}
\definecolor{labelkey}{rgb}{1,0,0}

\usepackage{cancel}

\usepackage{tikz}

  \mathchardef\ordinarycolon\mathcode`\:
  \mathcode`\:=\string"8000
  \begingroup \catcode`\:=\active
    \gdef:{\mathrel{\mathop\ordinarycolon}}
  \endgroup

\theoremstyle{plain}
\newtheorem{thm}{Theorem}[section]
\newtheorem{lem}[thm]{Lemma}

\newtheorem{cor}[thm]{Corollary}

\newcommand\eps{\varepsilon}
\newcommand\mi{\mathrm{i}}
\renewcommand\rho{\varrho}
\renewcommand\phi{\varphi}
\renewcommand\Re{\mathrm{Re}}
\renewcommand\Im{\mathrm{Im}}

\begin{document}

\title{\LARGE{Real chaos and complex time}}
\vspace{1cm}
{\subtitle{}
	\vspace{1ex}
	{}}\vspace{1ex}

\author{
 \\
\emph{-- Dedicated to the dear memory of Claudia Wulff  --}\\
{~}\\
Bernold Fiedler*\\
\vspace{2cm}}

\date{\small{version of \today}}
\maketitle
\thispagestyle{empty}

\vfill

*\\
Institut für Mathematik\\
Freie Universität Berlin\\
Arnimallee 3\\ 
14195 Berlin, Germany


\newpage
\pagestyle{plain}
\pagenumbering{roman}
\setcounter{page}{1}

\begin{abstract}
\noindent
Real vector fields $\dot{z} = f(z)$ in $\mathbb{R}^N$ extend to $\mathbb{C}^N$, for complex entire $f$.
We do not impose restrictions on the dimension $N$.
The homoclinic orbit $z(t)=\Gamma(t) :=1-3/ \cosh^2(t/\sqrt{2})$ to the equilibrium $z=1$ of the pendulum $\ddot{z}=z^2-1$ is a planar example.
Note the double poles of $\Gamma(t)$ at complex times $t/\sqrt{2}= \mi  (k+\tfrac{1}{2}) \pi$, for integer $k$.

\smallskip\noindent
Complex entire $f$ imply exponentially small upper bounds 
\begin{equation*}
\label{*}
C_\eta \exp(-\eta/\varepsilon) \tag{*}
\end{equation*}
on homoclinic splittings under discretizations of step size $\varepsilon>0$, or under rapid forcings with that time period.
Here $|\Im\, t|\leq \eta$ describes any complex horizontal strip where the complex time extension of $\Gamma(t)$ is analytic.
The phenomenon relates to adiabatic elimination, infinite order averaging, invisible chaos, and backward error analysis.
Claudia Wulff has studied symplectic variants.

\smallskip\noindent
However, what if $\Gamma(t)$ itself were complex entire?
Then $\eta$ could be chosen arbitrarily large.
We consider homoclinic or heteroclinic orbits $\Gamma(t)$ which connect limiting hyperbolic equilibria $f(v_\pm)=0$, for real $t\rightarrow\pm\infty$. 
For the diagonalizable linearizations $f'(v_\pm)$, we assume all unstable eigenvalues at $v_-$\,, and all stable eigenvalues at $v_+$\,, to be real. In addition we assume their nonresonance, separately at $v_\pm$\,.
We then show the existence of singularities of $\Gamma(t)$ in complex time $t$.
In that sense, real connecting orbits are accompanied by finite time blow-up, in imaginary time.
Moreover, the singularities bound admissible $\eta$ in exponential estimates \eqref{*}. 

\smallskip\noindent
The cases of complex or resonant eigenvalues are completely open.
We therefore offer a \textbf{1,000\,\euro\ reward} to any mathematician, up to and including non-permanent PostDoc level, who first comes up with a complex entire homoclinic orbit $\Gamma(t)$, in the above setting. 
Such an example would exhibit ultra-exponentially small separatrix splittings, and ultra-invisible chaos, under discretization.

\smallskip\noindent
The family $z(t)=\gamma(t)$ of periodic orbits generated by non-entire homoclinic orbits may well contain complex entire periodic orbits.
We provide a time reversible example.
In such cases, we predict ultra-sharp Arnold tongues, alias ultra-invisible phaselocking, under discretization.

\end{abstract}

\newpage
\tableofcontents


\newpage
\pagenumbering{arabic}
\setcounter{page}{1}

\section{Introduction} \label{Int}

\numberwithin{equation}{section}
\numberwithin{figure}{section}
\numberwithin{table}{section}

\subsection{Flows and iterations}\label{Flow}
The discrepancies between group actions $t\mapsto\Phi^t$ of diffeomorphisms on a real phase space $X$, for continuous real time $t\in\mathbb{R}$ versus discrete time $t\in \mathbb{Z}\eps$, are of fundamental interest in differentiable dynamical systems.
One-step numerical discretizations, for example, are the attempt to approximate autonomous flows, in continuous time, by cleverly chosen diffeomorphisms, for small time steps $\eps>0$.
Iteration typically addresses approximation over finite time intervals.
Averaging of non-autonomous systems, forced with small or finite time period $\eps$, on the other hand, attempts to approximate the associated stroboscopic period maps by interpolating autonomous flows.

By definition, such group actions satisfy 
\begin{equation}
\label{flow}
\Phi^t\circ\Phi^s= \Phi^{t+s}, \qquad \Phi^0=\mathrm{id},
\end{equation}
for all $t,s$.
In discrete time, integer multiples $t=n\eps$ just amount to $n$-th iterations $\Phi^{n\eps} = (\Psi^\eps)^n$ of some fundamental map $\Psi^\eps$.
Henceforth, in this notation, we will reserve $\Phi^t$ to denote differentiable flows in continuous time $t$, real for now, and complex later.

Standard existence and uniqueness theorems then identify differentiable flows with solutions $x=x(t)=\Phi^t(x_0)$ of initial value problems for autonomous \textbf{O}rdinary \textbf{D}ifferential \textbf{E}quations (ODEs) 
\begin{equation}
\label{ode}
\begin{aligned}
\dot{x}&=f(x),\\
x(0)&=x_0
\end{aligned}
\end{equation}
on $x\in X$, at least locally in time.
Here the autonomous, i.e. time-independent, vector field $f(x):= \partial_t\Phi^t(x)$ is evaluated at $t=0$. 
For further details see \cite{Hartman,Arnoldode,ArnoldODE, HirschSmaleDevaney, Devaney} and many other standard ODE references.
To keep notation simple, we only address globally defined flows, for now.
This excludes blow-up of solutions in finite time, which we will encounter later.
Adaptations to mere local flows are straightforward.

\subsection{Time discretization}\label{Dis}
For small discrete $\eps>0$, on the other hand, one-step numerical discretizations of a flow $\Phi^t$ define approximations 
\begin{equation}
\label{dis}
\Psi^\eps \ := \ \Phi^\eps + \mathcal{O}(\eps^{p+1}).
\end{equation}
The Landau symbol $\mathcal{O}$ indicates locally uniform boundedness after division by the argument, here $\eps^{p+1}$.
The positive integer $p$ is called the \emph{order} of the discretization scheme.
In other words, the Taylor expansions of the discretization $\Psi^\eps$, and of the time-$\eps$ stroboscope map $\Phi^\eps$, with respect to $\eps$, coincide up to order $p$ at $\eps=0$. The explicit Euler scheme
\begin{equation}
\label{euler}
\Psi^\eps(x) \ := \ x+\eps f(x),
\end{equation}
for example, has $p=1$.
The ultimate no-error ``discretization'' $\Psi^\eps \ = \ \Phi^\eps$, by the true, original flow $\Phi^t$ itself, constitutes any numerical analyst's dream.
Alas, it rarely happens.
See any good reference on ODE numerics, like the excellent monograph \cite{HairerLubichWanner}. 
For reduction of multi-step and variable-step methods to fixed single steps, see also \cite{Kirchgraber-multi, Stoffervar, Stoffer-multi}.

\emph{Backward error analysis} takes a subtly different perspective; see for example \cite{Hairer, HairerLubich, Reich, Matthiesback, HairerLubichWanner, CohenHairerLubich}.
It seeks an \emph{interpolating ODE flow} $\widetilde\Phi^t$, for given discretization $\Psi^\eps$ of $\Phi^\eps$, such that the approximation \eqref{dis} holds with much, much higher precision, for the ``wrong'' flow $\widetilde\Phi^t$ instead of the true, original flow $\Phi^t$.
For real analytic vector fields $f$, typical estimates are \emph{exponential}, i.e.
\begin{equation}
\label{exp}
\Psi^\eps \ = \ \widetilde\Phi^\eps + \mathcal{O}(\exp(-\eta/\eps)),
\end{equation}
locally uniformly for some fixed $\eta>0$, and $\eps\rightarrow 0$.
Such exponential smallness estimates are sometimes also called \emph{``beyond finite order''} or \emph{``of infinite order''} in $\eps$.

Let $\tilde{f}$ denote the autonomous ODE vector field generating the ``wrong'' flow $\widetilde{\Phi}^t$.
Then the study of numerical approximations decomposes into two steps.

First, there is a discrepancy $\mathcal{O}(\eps^{p})$ between the original ODE vector field $f$ and $\tilde{f}$; see \eqref{dis}, \eqref{exp}.
On the bright side, at least the autonomous flow character is preserved exactly.
In cases of structural stability of $f$, for example, $C^0$ orbit-conjugacy of the two ODEs holds, for all small $\eps$. This includes hyperbolic flows, e.g., with chaotic dynamics \cite{PalisdeMelo,Devaney, HasselblattKatok}, and uniformly addresses an infinite time horizon.

Second, there only remains an exponentially small discrepancy between the flow map $\widetilde{\Phi}^\eps$, of the ``wrong'' ODE $\tilde{f}$, and the $p$-th order discretization $\Psi^\eps$, of the ``true'' ODE map $\Phi^\eps$.
In other words, that ``wrong'' ODE is well-approximated by the original numerical discretization, over exponentially long time horizons $t= \mathcal{O}(\eta/\eps)$.
Combined with structural stability techniques like orbit conjugacy, this makes for extremely precise simulations, uniformly over exponentially large times, under quite sloppy discretizations.
Symplectic discretization, for example, well preserves the Hamiltonian $\widetilde{H}$ of the ``wrong'' flow $\widetilde{\Phi}^t$, although incorrect of order $\eps^p$ as a function, over exponentially long times.

\subsection{Rapid forcing and normal forms}\label{Nor}
In the Hamiltonian averaging setting of rapid periodic forcings, exponential ODE estimates \eqref{exp} have first been established by Neishtadt under the name of \emph{adiabatic elimination} \cite{Neishtadt}, also called averaging \emph{beyond finite order} or \emph{of infinite order}.
That setting is equivalent to symplectic discretization, by Hamiltonian symplectic interpolation of $\Psi^\eps$; see for example \cite{MoserSymplectic, MoserPhysics}.
For further numerical developments see \cite{HairerLubich}. 
PDE results for analytic semigroups with compact resolvents have first  been obtained by Karsten Matthies \cite{MatthiesDiss}.
Under high spatial regularity assumptions of Gevrey class, he obtained splitting exponents $-\eta\eps^{-1/3}$. 
See \cite{MatthiesScheel} for generalizations to Hamiltonian PDEs like hyperbolic wave equations,
\cite{Matthiesell} for exponential homogenization of elliptic PDEs in cylinder domains,
and \cite{Matthiesquasiperiodic} for PDEs with rapid quasi-periodic forcings.
With Marcel Oliver, Claudia Wulff has extended exponential backward error analysis to infinite-dimensional symplectic settings \cite{OliverWulff, WulffOliver}.
In a somewhat related spirit, she has addressed almost invariant symplectic slow manifolds under singular perturbations, up to exponential precision in the singular parameter $\eps$; see \cite{KristiansenWulff} and compare \cite{IoossLombardi}.
From a conceptually much broader perspective, she also worked on symmetry breaking of relative periodic orbits, up to her last paper \cite{ChillingworthLauterbachWulff}.

A slightly different flavor of exponential smallness arises from normal form calculations involving equilibria with purely imaginary eigenvalues.
See for example the monograph \cite{Lombardi}, with beautiful applications to the dynamics of water waves.
Formal normal forms first provide an $\mathbb{S}^1$-equivariance, to any finite algebraic order.
Beyond finite order, deviations from angular $\mathbb{S}^1$- equivariance may then be thought of as a rapidly oscillating perturbation, in the angular phase direction.
On large homoclinic orbits, this induces a rapid forcing, and generates large orbits which are transversely homoclinic to small periodic solutions.
Although temptingly similar, in flavor, we cannot pursue the intriguing aspect of perturbations beyond finite order normal form further, in the present paper.

\subsection{Lie algebras, Lie groups, and exponentiation}\label{Lie}
From an abstract and formal point of view, we may consider ODE vector fields $f$ as the Lie algebra, i.e. the infinitesimal version, of the local Lie group of near-identity diffeomorphisms $\Psi^\eps$.
Exponentiation becomes the passage from $f$ to the flow $\Phi^t$.
Discretization encounters non-surjectivity of exponentiation.

As a trivial first example, let us consider linear ODE flows $\dot{x}=Ax$, and the discretization of their flow maps $\Phi^\eps=\exp(\eps A)$ by near-identity linear diffeomorphisms $\Psi^\eps=\mathrm{id}+\eps B$. 
For small $\eps$, there need not be any discrepancy \eqref{exp} at all, be it quantitative or just qualitative.
Locally, indeed, we may always invert exactly,
\begin{equation}
\label{LieAlgebra}
\begin{aligned}
\Psi^\eps &= \mathrm{id}+\eps B =  \exp(\eps \tilde{A}) \quad \Longleftrightarrow\\
\tilde{A} &= \eps^{-1} \log (\mathrm{id}+\eps B) = B+\ldots\,.
\end{aligned}
\end{equation}
In other words, for any linear discretization $\mathrm{id}+\eps B$ sufficiently near identity, there exists a ``wrong'' ODE $\dot{x}=\tilde{A}x$ which is integrated exactly, without error, by that discretization.
In general, that ODE will differ from the true, original ODE $\dot{x}=Ax$, which gave rise to the discretization in the first place.
Standard Floquet theory cautions us that exact inversion may already fail for real matrices and finite times like $\eps=1$.

In general, exponentiation provides a local diffeomorphism between any finite-dimensional Lie algebra and its Lie group.
For example, we may truncate the relation between an ODE vector field $f$ and its flow $\Phi^t$ to some arbitrarily large, but finite, polynomial degree $k$. Abstractly, the relation then amounts to exponentiation: from the finite-dimensional Lie algebra of vector fields, truncated at order $k$, into the Lie group of truncated near-identity local diffeomorphisms, under (truncated) composition \eqref{flow}.
This already suggests that approximations \eqref{exp} might be possible, beyond any finite order $k$.

\subsection{Roots of diffeomorphisms}\label{Root}
Without truncation, however, fundamental discrepancies appear between the dynamics of iterations $n\mapsto (\Psi^\eps)^n$, in discrete time $n\in\mathbb{Z}$, and the stroboscopic view $n\mapsto (\Phi^\eps)^n=\Phi^t$ of the original flow, evaluated at discrete times $t=n\eps$.

Consider any diffeomorphism $\Psi$, for example.
We call $\Psi$ \emph{infinitely divisible}, if we can define $n$-th roots $\Psi^{1/n}$ of $\Psi$ in the diffeomorphism group, such that 
\begin{equation}
\label{root}
(\Psi^{1/n})^n=\Psi,
\end{equation}
for any integer $n$.
The time-1 diffeomorphism $\Psi=\Phi^1$ of any flow $\Phi^t$, for example, is infinitely divisible with $\Psi^{1/n}:= \Phi^{1/n}$.

Conversely, \eqref{root} also defines rational exponents $\Psi^{m/n}$.
Assuming continuity, we obtain an associated flow $\widetilde{\Phi}^t := \Psi^t$ for all real $t$.
Assuming differentiability, we can define an associated autonomous ODE $\dot{x}=\tilde{f}(x)$ with $\Psi=\widetilde{\Phi}^1$ as its time-1 stroboscope map.

The simplest logistic equation $\dot{x}=x(1-x)$ on $x\in X=[0,1]$ provides an interesting nonlinear example.
Already the beautifully detailed analysis by Belitskii and Tkachenko has addressed  diffeomorphisms $\Psi$ on $X=[0,1]$ with $x=0,1$ as the only fixed points.
Both fixed points are assumed hyperbolic, e.g. $\Psi'(0)>1>\Psi'(1)>0$.
Generically, it turns out, such diffeomorphisms do not possess any root of order $n>1$.
More precisely, section 2.6 of \cite{BelitskiiTkachenko} determines the moduli space of conjugacy classes of such $\Psi$ in  $C^k$, via certain periodic matching functions on a fundamental interval like $I=[\tfrac{1}{2}, \Psi(\tfrac{1}{2})]$.
The matching functions keep track of the mismatch between the two locally linearizing diffeomorphisms, at the hyperbolic fixed points $x=0$ and $x=1$, respectively, after propagation to the overlap interval $I$.
 A root $\Psi^{1/n}$ exists if, and only if, the matching function cycles through  $n$ (not necessarily minimal) periods over the fundamental interval $I$.
Time-1 maps of flows, in particular, are therefore characterized by constant matching functions.
 
The absence of $n$-th roots for generic diffeomorphisms $\Psi$ on general compact manifolds is related to Smale's \emph{centralizer conjecture}.
Consider the group $\mathcal{G}$ of iterates $\Psi^n,\ n\in\mathbb{Z}$.
The \emph{centralizer} $\mathcal{Z}=\mathcal{Z}(\mathcal{G})$, as usual, consists of all diffeomorphisms $\Phi$ which commute with each element of $\mathcal{G}$.
Because $\mathcal{G}$ is generated by $\Psi$, in our case, the centralizer coincides with the commutator of $\Psi$, i.e. with all $\Phi$ such that $\Phi\circ\Psi=\Psi\circ\Phi$.
In particular $\mathcal{Z}\supseteq\mathcal{G}$.

As dynamical systems problem $12$ in his 1998 list, Smale has conjectured $\mathcal{Z}=\mathcal{G}$, i.e. triviality of the centralizer, for generic $\Psi$.
See \cite{SmaleProblems}, dating back to \cite{SmaleCentralizer}.
In the $C^1$-topology, the conjecture has been proved in \cite{SmaleRoot}.

Generic absence of $n$-th roots is a straightforward corollary.
First note $\Phi\in\mathcal{Z}$, in case any $n$-th root $\Phi:=\Psi^{1/n}$ of $\Psi$ did exist.
Indeed $\Phi\circ\Psi=\Phi^{n+1}=\Psi\circ\Phi$.
For generic $\Psi$, moreover, $\mathcal{G}$ is not finite cyclic, due to the Kupka-Smale theorem \cite{AbrahamRobbin,PalisdeMelo}.
In particular $\Phi\in\mathcal{Z}\setminus\mathcal{G}\neq\emptyset$, because $\Psi$ is not idempotent.
This contradiction to Smale's centralizer conjecture proves generic absence of $n$-th roots.
In the $C^0$-topology, in contrast, structural stability may prevail -- with or without roots.
See for example the elementary discussion of hyperbolic torus diffeomorphism and their roots in \cite{FiedlerAnosov}. 

For a precursor on centralizers of circle diffeomorphisms in the $C^k$-topology, $k\geq 2$, see also \cite{KopellRoot}.
Circle diffeomorphisms $\Psi\in C^2$ on $X=\mathbb{S}^1=\mathbb{R}/\mathbb{Z}$, however, exhibit a well-known peculiarity.
For irrational rotation numbers $\rho$, they are $C^0$-conjugate to translations $x\mapsto x+\rho \in \mathbb{S}^1$, i.e. to the time-1 map of the translational flow $\dot{x}=\rho$.
In particular, they are $C^0$ infinitely divisible.
For rational rotation numbers $\rho=r/s$, however, $\Psi$ possesses periodic orbits of minimal period $s$.
If these are all hyperbolic, then the Belitskii-Tkachenko obstruction to divisibility applies to any interval of adjacent fixed points of $\Psi^s$.
Moreover, the rotation number remains pinned at the plateau value $\rho=r/s$, under $C^1$-small perturbations.
For circle maps $\Psi$, this phenomenon has been popularized under the name \emph{devil's staircase} or \emph{Cantor function}, in generic one-parameter families $\rho(\lambda)$ \cite{Brunovsky, Aubry, Devaney}, and in music \cite{Ligeti}.
In generic two-parameter families, it arises in the guise of  \emph{Arnold tongues} \cite{ArnoldODE, Kuznetsov, Devaney}.
Similar phenomena occur for discretizations $\Psi^\eps$ of hyperbolic periodic orbits $\Phi^T(x)=x$ with minimal period $T>0$.
See section \ref{UsA}.

\subsection{Transverse homoclinic orbits}\label{Trv}
Stable and unstable manifolds have coincided in the above one-dimensional examples of heteroclinic orbits $\Psi^n(x)$ between adjacent hyperbolic fixed points of adjacent unstable dimensions.
Heteroclinic orbits just filled the in-between interval, monotonically.
Phase spaces $X$ of dimension $N\geq 2$ admit much more interesting dynamics.
For example, we may consider nonstationary orbits
\begin{equation}
\label{homhet}
x(t)=\Gamma(t) \rightarrow v^\pm, \quad \textrm{ for } t\rightarrow\pm\infty
\end{equation}
between hyperbolic equilibria or fixed points $x=v^\pm$.
We call $\Gamma$ \emph{homoclinic}, in case $v^+=v^-=v$\,, and \emph{heteroclinic} otherwise.

For ODE flows, and in absence of further structure like a first integral, symplecticity, or reversibility, homoclinicity typically requires the adjustment of an additional scalar parameter $\lambda$, i.e.
\begin{equation}
\label{odef}
\dot{x}=f(\lambda,x).
\end{equation}
Assume $f(\lambda,v)=0$, for all $\lambda$.
The scalar real parameter $\lambda$ is introduced to ensure that stable and unstable manifolds of the hyperbolic equilibrium $v$ cross each other at nonvanishing speed, as $\lambda$ increases through $\lambda=0$.
Indeed, the stable and unstable manifolds of $v$ still coincide, along the homoclinic curve $\Gamma$ and, therefore, \emph{cannot be transverse}.
The same remark applies to heteroclinic orbits, in case the unstable dimensions of the hyperbolic equilibria $v^\pm$ coincide.
For generic $C^k$-diffeomorphisms $\Psi$, in contrast, homoclinic (and heteroclinic) intersections \emph{are transverse}.
Together with generic hyperbolicity of fixed points and periodic orbits, this constitutes the celebrated Kupka-Smale theorem; see for example \cite{AbrahamRobbin, PalisdeMelo}.

Transverse homoclinic orbits persist under $C^1$-small perturbations.
By the celebrated Smale horseshoe construction, they also give rise to nearby shift dynamics; see for example \cite{Moser, GuckenheimerHolmes, Palmer, PalisTakens, Devaney} and the references there.
For a delicate $C^0$-invariant associated to structurally stable transverse homoclinicity in Anosov diffeomorphisms of 2-tori, see also \cite{FiedlerAnosov}.
Homoclinic transversality versus nontransversality, of course, is yet another fundamental qualitative discrepancy between dynamics in continuous and discrete time.
In the Lie group language of section \ref{Lie}, this discrepancy is yet another indicator of local non-surjectivity of the exponential map from vector fields to near-identity diffeomorphisms.

\begin{figure}[t]
\centering \includegraphics[width=0.6\textwidth]{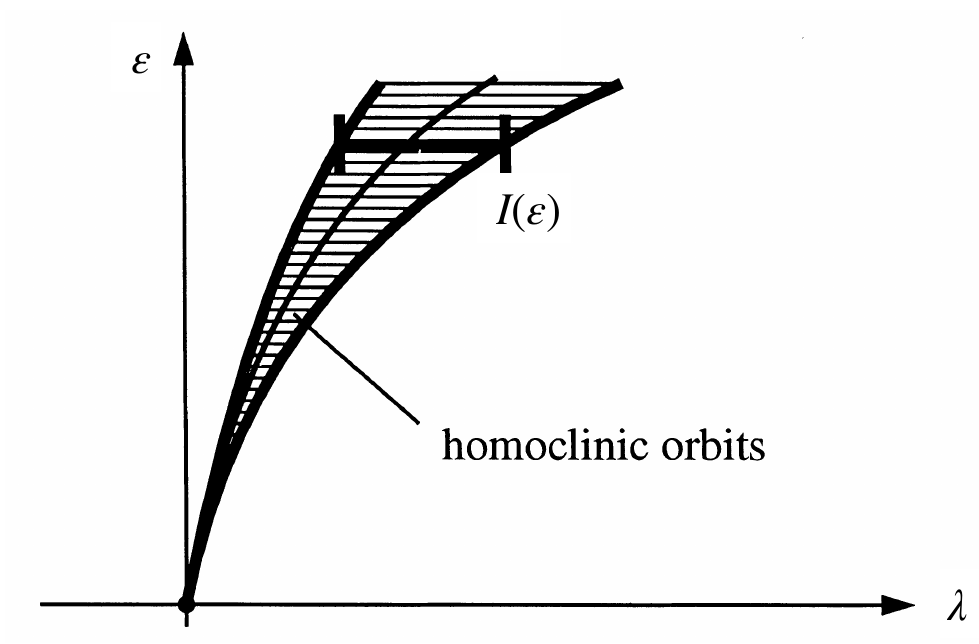}
\caption{\emph{
Schematic splitting region for a homoclinic curve $\Gamma$ of the ODE-flow \eqref{odef} (hashed), under discretization with step size $\eps$ or, equivalently, under $\eps$-periodic non-autonomous forcing \eqref{odeg}; see \cite{FiedlerScheurle}.
At fixed levels of $\eps$, the horizontal splitting intervals $\lambda\in I(\eps)$ mark parameters $\lambda$ for which single-round homoclinic orbits occur near $\Gamma$.
In the analytic case, backward error analysis indicates exponential smallness \eqref{split} of the splitting width $\ell(\eps)$ of $I(\eps)$, for some small $\eta>0$.
Established splitting exponents $\eta$ are in fact strictly bounded by the half-width of the analyticity strip of the complex extension of $\Gamma(t)$, in the imaginary time direction.
Exponential smallness estimates apply, analogously, to the splitting region of heteroclinic orbits between hyperbolic equilibria $v^\pm$ of equal unstable dimension.
For an interpretation of the same figure in terms of resonant discretization of periodic orbits, and ultra-invisible Arnold tongues, see also theorem \ref{thmultra-exp} and fig. \ref{fig3}.
For better visibility, the horizontal width $\ell(\eps)$ of the exponentially flat splitting region has been exaggerated, in our schematic illustration.\\
In the homoclinic case, we call the recurrent dynamics in the resulting chaotic region \emph{``invisible chaos''}.
For complex entire homoclinic orbits $\Gamma(t)$, the exponent $\eta$ could be chosen arbitrarily large.
Such ultra-exponentially small splitting, although currently elusive, would lead to \emph{``ultra-invisible chaos''}.
}}
\label{fig1}
\end{figure}

\begin{center}
But: how about \emph{quantitative} discrepancies?
\end{center}

Backward error analysis \eqref{dis}, \eqref{exp} suggests \emph{splitting intervals} $\lambda\in I(\eps)$ \emph{of exponentially small length} $\ell(\eps)$ where homoclinic transversality may prevail, i.e. an upper estimate
\begin{equation}
\label{split}
\ell(\eps) \leq C_\eta \exp(-\eta/\eps),
\end{equation}
for some fixed constants $\eta, C_\eta>0$ and all small $\eps \searrow 0$.
See fig.~\ref{fig1}.
To be more precise, we consider \emph{1-homoclinic} here, alias \emph{single-round homoclinic orbits}, which remain uniformly close to the original homoclinic orbit $\Gamma(t)$, in real time.

\subsection{Analytic extension and rapid forcing}\label{Rap}
To specify the admissible range of exponents $\eta>0$ in \eqref{split}, more quantitatively, we are led to consider extensions of ODE \eqref{odef} to complex time.
For a primer on several complex variables, see appendix A of \cite{Ilyashenko} and the references there.
Just recall that a complex-valued function of one or several complex variables is called \emph{analytic} in an open complex domain $\mathcal{D}\subseteq\mathbb{C}^N$, if it is represented by its local Taylor series there.
Equivalently, the function is \emph{holomorphic} in $\mathcal{D}$, i.e. complex differentiable in the sense of the Cauchy-Riemann equations.
In case $\mathcal{D}=\mathbb{C}^N$, the function is called \emph{complex entire}.
For disambiguation, we avoid PDE terminology here, which often calls global or globally bounded real solutions $x(t),\ t\in\mathbb{R}$, ``entire'' or ``eternal''.

Under local analyticity assumptions, in fact, exponentially small splitting estimates \eqref{split} hold for \emph{any} fixed $\eta>0$ such that the original homoclinic orbit $\Gamma(t)$ extends, complex analytically, to a neighborhood of the closed complex horizontal strip
\begin{equation}
\label{horstrip}
|\Im\, t| \leq \eta\,.
\end{equation}
Here and below, $\sigma=\Re\,t$ and $\tau=\Im\,t$ denote the real and imaginary parts of complex time $t$.

In \cite{FiedlerScheurle}, the estimate \eqref{split}, \eqref{horstrip} has actually been proved for homoclinic splittings in a setting of rapid periodic forcings
\begin{equation}
\label{odeg}
\dot{x}(\sigma)=f(\lambda,x(\sigma)) + \eps^p g(\lambda,\eps,\sigma/\eps,x(\sigma))\,.
\end{equation}
The functions $f,g$ are assumed to be complex entire in $\lambda, x,\eps$, and real for real arguments. 
In the variable $\beta=\sigma/\eps$ of the rapid forcing $g=g(\lambda,\eps,\beta,x)$, we only assume smoothness and period 1.
The origin $x=0$ is assumed to be a hyperbolic equilibrium of \eqref{odef} with nondegenerate homoclinic orbit $\Gamma(\sigma)$, as described above.
Under a spectral nondegeneracy condition, and in the planar case, this leads to \emph{invisible chaos}: the chaotic region accompanying transverse homoclinics is exponentially thin, both, in parameter space and along the homoclinic loop.
Equivalence of time $\sigma=\eps$ stroboscope maps $\Psi^\eps$ of \eqref{odeg} with arbitrary analytic $p$-th order discretizations $\Psi^\eps$ of the autonomous flow of \eqref{odef} has been established.
Therefore, our results also apply to discretizations of \eqref{odef}.
See \cite{FiedlerScheurle} for complete details, and \cite{Matthieshom} for generalizations of exponentially small homoclinic splitting estimates to analytic PDE semigroups.

\subsection{Complex entire homoclinic orbits -- a 1,000\,\euro\ question}\label{1000}
At $\eps=0$, in particular, we require analyticity of the homoclinic orbit $\Gamma$ for times $t=\sigma+\mi \tau$ in complex horizontal strips \eqref{horstrip}.
This raises the question what would happen to the upper exponential estimates \eqref{split} for \emph{complex entire homoclinic orbits} $\Gamma(t)$, i.e. under global analyticity for all complex $t\in\mathbb{C}$.
Then exponential estimates \eqref{split} hold, simultaneously for \emph{all} $\eta>0$, sufficiently small $\eps>0$, and suitable constants $C_\eta$\,.
Taking logarithms, we obtain $-\log \ell(\eps) \geq c_\eta +\eta/\eps$, with real $c_\eta:= -\log C_\eta$\,.
In other words, we obtain upper splitting estimates
\begin{equation}
\label{ultrasplit}
-\log \ell(\eps) \geq c(1/\eps) >0 ,
\end{equation}
for a suitable convex function c of unbounded positive slope.
We call this phenomenon \emph{ultra-exponentially small splitting of separatrices}.
In the planar setting of \cite{FiedlerScheurle}, section 5, such splitting would be accompanied by \emph{ultra-invisible chaos}.

Alas, can it actually happen?
The author personally offers a
\begin{center}
\textbf{1,000\,\euro\ reward}
\end{center}
for settling this question.

More precisely, a \emph{positive answer} would require any explicit example of a real nonstationary homoclinic orbit $\Gamma(t)$ to a hyperbolic equilibrium, for a complex entire ODE \eqref{ode} on $X=\mathbb{C}^N$ with $f$ real for real arguments, such that  $\Gamma(t)$ is complex entire for $t\in\mathbb{C}$. 
Such an example would initiate many challenges, e.g., for numerical exploration.
Eventually, it might lead towards a whole new theory to address ultra-exponentially small splitting behavior under discretization.

A \emph{negative answer} would prove that such complex entire homoclinic orbits $\Gamma(t)$ cannot exist.
In the planar case $N=2$, indeed, the answer is known to be negative; see the discussion section in \cite{FiedlerScheurle}. 
See also earlier results on one-dimensional unstable manifolds in \cite{Ushiki-2, Ushiki-N}.
In the planar Hamiltonian case, hyperbolicity of the equilibrium can also be dropped.
For the special case of scalar second order pendulum equations 
\begin{equation}
\label{pendulum}
\ddot z + \mathbf{q}(z) = 0
\end{equation}
with nonlinear complex entire $\mathbf{q}$, already Rellich \cite{Rellich} in 1940 has most elegantly noticed the complete absence of \emph{any} non-constant complex entire solutions $z(t)$; see also \cite{Wittich2}.
Wittich has generalized this result to scalar nonautonomous $m$-th order equations which are a polynomial $\mathbf{p}$ in $t,z,\dot z,\ldots,z^{(m)}$ plus a non-polynomial complex entire nonlinearity $\mathbf{q}(z)$,
\begin{equation}
\label{pendulumm}
\mathbf{p}(z^{(m)},\ldots,z,t)+\mathbf{q}(z)=0;
\end{equation}
see \cite{Wittichm}.
Under an elementary extra condition on $\mathbf{p}$, his result includes the purely polynomial case $\mathbf{q}=0$.
The only complex entire solutions $z(t)$ are then polynomial in $t$, and hence neither homo- or heteroclinic, nor non-constant periodic.
In particular, the Wittich condition excludes harmonic solutions like $\dot z^2+z^2=1$ or $\ddot z+z=0$.
The Wittich condition is violated in our autonomous, time reversible, second order example \eqref{rev} below -- which promptly features a periodic cosine solution \eqref{revper}.

In order to protect well-established colleagues against potentially distracting personal financial interests, the prize will be awarded to the first solution by anyone \emph{without} a permanent position in academia.
Priority is defined by submission time stamp at arxiv.org or equivalent repositories.
Subsequent confirmation by regular refereed publication is required.

\subsection{Real nonresonant spectra}\label{Res}
The goal of our present paper is a \emph{negative answer} in the case where stable and unstable spectra of the hyperbolic equilibrium $x=v$ are real and \emph{nonresonant, separately}.
Under the above analyticity assumptions on $f$ in \eqref{ode} , let $A^\pm$ denote the block decomposition of the linearization $A=f'(v)$ at $f(v)=0$, such that $A^-$ represents all $N^-$ stable eigenvalues $\mu$ of $A$, and  $A^+$ represents all remaining $N^+=N-N^- $ unstable eigenvalues of $A$.
See also \eqref{odeclin}, \eqref{odecpm}.
In other words, 
\begin{equation}
\label{specpm}
\pm\,\Re\,\mu >0, \textrm{ on } A^\pm\,.
\end{equation}
We assume $A^\pm$ to be diagonalizable, i.e. algebraic and geometric spectral multiplicities coincide.
\emph{Stable nonresonance} requires that all eigenvalues $\mu^-_j\in \mathrm{spec}\,A^-$ satisfy the nonresonance conditions
\begin{equation}
\label{nonres-}
\mu^-_{j_0} \neq \sum_j n_j \,\mu^-_j\,,
\end{equation} 
for any index $j_0$, and any multi-index $\mathbf{n}=(n_1,\ldots,n_{N^-})\in\mathbb{N}_0^{N^-}$ of order $|\mathbf{n}|:=n_1+\ldots+n_{N^-}\geq 2$.
Analogously, \emph{unstable nonresonance} on $\mathrm{spec}\,A^+$ requires 
\begin{equation}
\label{nonres+}
\mu^+_{j_0} \neq \sum_j n_j \,\mu^+_j\,.
\end{equation}

\subsection{Real nonresonant homoclinic and heteroclinic orbits are not complex entire}\label{Thms}
We are now ready for our main result: the exclusion of complex entire homoclinic orbits $\Gamma(t)\rightarrow v$, for $\Re\,t\rightarrow\pm\infty$, in the case of real, separately nonresonant eigenvalues.

\begin{thm}\label{thmhom}
In the above setting, assume semisimple, real, and separately nonresonant eigenvalues $\mu_j^\pm$ of the diagonalizable linearizations $A^\pm$, at the hyperbolic equilibrium $v$; see \eqref{nonres-}, \eqref{nonres+}.\\
Then complex entire homoclinic orbits $\Gamma(t)$ to $v$ cannot exist.
\end{thm}

Let us take a moment to critically assess the relevance of the theorem for the homoclinic splitting problem.
In dimension $N=2$, eigenvalues are automatically real and nonresonant.
This recovers the planar results of \cite{Ushiki-2, FiedlerScheurle}, in the hyperbolic case.
In higher dimensions, of course, our negative theorem \ref{thmhom} still permits complex entire homoclinic candidates under real spectral resonances, or involving conjugate complex pairs of eigenvalues.
The elimination of real (nonresonant) spectrum specifically points at Shilnikov homoclinic orbits of saddle-focus type, as potentially complex entire 1,000\ \euro\ candidates, and also at homoclinic orbits to saddle-node equilibria.
See for example \cite{Shilnikov}.

A more fundamental critique concerns the relevance of mere upper exponential estimates like \eqref{split}.
Such one-sided estimates do not quantify the actual separatrix splitting, unless complemented by \emph{lower estimates} of exponential type.
For complex entire planar diffeomorphisms, the presence of nonzero splittings is (rather) implicit in \cite{Ushiki-2}.
For the rapidly periodically forced pendulum, a lower exponential estimate has been established by \cite{Scheurleetal} in 1991; see also \cite{Fontich} and the references there.

In many planar systems, mostly of symplectic type, a longstanding effort by Lazutkin and coworkers has culminated in a large body of admirable and diligent work by Gelfreich; see \cite{Lazutkin} and the surveys \cite{Gelfreich01, Gelfreich02}.
That very fine analysis has established \emph{asymptotic expansions}, for separatrix splittings, of the form $c(\epsilon) \exp (-\eta/\epsilon)$ with a formally meromorphic coefficient $c(\epsilon)$.
Here $\eta$ stands for the exact distance of the complex poles of $\Gamma(t)$ from the real axis, in the complex time plane $t\in\mathbb{C}$.
The expansion parameter $\epsilon$ is usually closely and explicitly related to $\eps$ itself.
Among other questions, these results settled many exponential splitting problems which had first been encountered by Poincaré's famous error, and discovery, in the 3-body problem of celestial mechanics \cite{Poincare3body}, and which had remained open for a century.

For certain quasi-periodic rapid forcing functions $g$ and exponential estimates in fractional powers of $1/\eps$, see \cite{Delshams, Delshamswhiskered} and the references there.

Analogously to theorem \ref{thmhom}, we can exclude complex entire real heteroclinic orbits $\Gamma(t)\rightarrow v^\pm$ between distinct hyperbolic equilibria $v^\pm$, for $\Re\,t\rightarrow\pm\infty$. 
This time, we have to assume nonresonance \eqref{nonres-} for the stable part $A^-$ of the linearization $f'(v^+)$ at the target equilibrium $v^+$, and \eqref{nonres+} for the unstable part $A^+$ of the linearization $f'(v^-)$ at the source equilibrium $v^-$.

\begin{thm}\label{thmhet}
In the above modified heteroclinic setting, assume semisimple, real, and separately nonresonant eigenvalues $\mu_j^\pm$ of the diagonalizable linearizations $A^\pm$ at the hyperbolic equilibria $v^\mp$; see \eqref{nonres-}, \eqref{nonres+}.\\
Then complex entire heteroclinic orbits $\Gamma(t)$ from $v_-$ to $v^+$ cannot exist.
\end{thm}

As a corollary, we observe a curious connection between homo- or heteroclinicity of $\Gamma$, in real time $\sigma$, and \emph{blow-up in finite imaginary time} $\mi \tau$, for extensions to complex time $t=\sigma+\mi \tau$.
To our knowledge, such a remarkable link has first been discovered in a PDE context (!) by Stuke \cite{Stuke} in his brilliant dissertation thesis.
A much less explicit foreboding, in the different and very restrictive setting of complex entire planar diffeomorphisms, can be attributed to Ushiki \cite{Ushiki-2}.
For the special heteroclinic case $N^+=N^-=1$ in higher dimensions $N$, albeit of high codimension $N-1$, see also \cite{Ushiki-N}.

\begin{cor}\label{corblowup}
In the above settings, respectively, consider any real homoclinic or heteroclinic orbit $\Gamma(\sigma)$; see \eqref{homhet}. 
Assume semisimple, real, and separately nonresonant eigenvalues $\mu_j^\pm$ of the diagonalizable linearizations $A^\pm$, at real times $\sigma=\mp\infty$.\\
Then there exists some real homoclinic initial condition $\Gamma(\sigma_0)$, such that the complex analytic extension $\tau \mapsto \Gamma(\sigma_0+\mi \tau)$ blows up at finite imaginary times $\tau=\pm\tau^*\neq 0$.
\end{cor}

\subsection{Outlook on PDE blow-up}\label{PDE}
In corollary \ref{corblowup} we have viewed the absence of complex entire time extensions of real homoclinic or heteroclinic orbits $\Gamma(\sigma)$ as an indication of blow-up $|\Gamma(t)|\rightarrow\infty$ in finite complex time $t=\tau+\mi\sigma$.
In the sequel \cite{FiedlerStuke}, we further pursue such questions in a PDE context.
Following pioneering work by Masuda \cite{Masuda} in the spirit of Stuke \cite{Stuke}, we study real heteroclinic orbits $\Gamma(\sigma)$ from hyperbolic equilibria $v^-$ to $v^+$ in heat equations with a quadratic nonlinearity and Neumann boundary conditions.
Since the stable manifold at the target equilibrium $v^+$ is infinite-dimensional, Poincaré linearization is not available at $v^+$.
Assuming $v^+$ to be asymptotically stable, we still have to tackle unstable nonresonance at $v^-$.
For unstable dimensions of $v^-$ not exceeding 22, or within sufficiently fast unstable submanifolds of $v^-$, we still establish an analogue of blow-up corollary \eqref{corblowup}, at least for most quadratic nonlinearities.

\subsection{Overview}\label{Ove}
Our proofs of theorems \ref{thmhom} and \ref{thmhet} will be indirect.
Assuming the existence of a complex entire homo- or heteroclinic orbit $\Gamma(t),\ t\in\mathbb{C}$, under our standing spectral assumptions, we will reach a contradiction in section \ref{Pro}.
To build up, section \ref{Cot} provides some elementary background on flows in real versus complex time, including some examples, and our adaptation of local stable and unstable manifolds to the complex setting.
Section \ref{Liq} recalls and applies Poincaré linearization to these locally flow-invariant manifolds, under the spectral nonresonance conditions \eqref{nonres-}, \eqref{nonres+}.
The restriction to nonresonant real eigenvalues causes global quasi-periodicity of the flows on the complex local invariant manifolds, in the imaginary time direction.
Section \ref{Fos} therefore recalls Bohr's classical notion of almost-periodic Fourier series.
As a final technical tool, to reach the promised contradictions in section \ref{Pro}, we establish an invariance of the formal Fourier coefficients for any complex entire homo- or heteroclinic orbit $\Gamma(t)$, in section \ref{Foi}.
We conclude, in section \ref{UsA}, with a suggestion for an ultra-exponentially sharp Arnold tongue.

\subsection{Acknowledgment}\label{Ack}
This paper is dedicated to the dear memory of my doctoral daughter Claudia Wulff: to her uncompromising dedication to Science and Truth, to her moments of pure happiness and abandon, to her inspiration and charismatic liveliness against all odds.
I am also indebted to Vassili Gelfreich, for very patient explanations of his profound work on exponential splitting asymptotics, to Anatoly Neishtadt for enlightening conversations, to Jürgen Scheurle, and of course to Claudia herself, for their long-term interest.

\section{Complex time} \label{Cot}

In this section we discuss some basics concerning local flows and ODEs in complex time $t=\sigma+\mi\tau$.
See the first chapters of \cite{Ilyashenko} for a beautiful and complete introduction. 
We then review stable and unstable manifolds, in the complex time setting.

Our general setting, in the present paper, is based on autonomous ODEs \eqref{ode} with vector fields $f$ on $X=\mathbb{R}^N$. 
We assume $f$ to be given by a globally convergent real Taylor series.
In particular, the real ODE \eqref{ode} extends to a complexified ODE
\begin{equation}
\label{odec}
\begin{aligned}
\dot{z}&=f(z)\\
z(0)&=z_0 \in Z=\mathbb{C}^N.
\end{aligned}
\end{equation}
In other words, the complexified vector field $f:Z\rightarrow Z$ is complex entire.
Moreover $f$ commutes with complex conjugation, i.e. 
\begin{equation}
\label{fcc}
f(\bar{z})=\overline{f(z)}.
\end{equation}

We identify $Z=\mathbb{C}^N\cong \mathbb{R}^{2N}$ and abbreviate $z=x+\mi y=(x,y)$.
This provides a (local) flow $\Phi^t(z_0):=z(t)$ on $Z$, for real times $t$.
The local solutions are in fact real analytic, in $t$.
This extends \eqref{odec}, and the local flow $\Phi^t$, to complex times $t$.
In particular, the flow property \eqref{flow} extends to complex $s,t$, verbatim, at least locally.
Conversely, complex differentiation of $\Phi^t$ with respect to complex $t$ satisfies \eqref{odec}.

By restriction, \eqref{odec} defines standard real flows along any line $\sigma\mapsto t=e^{\mi \vartheta}\sigma,\ \sigma\in\mathbb{R}$, for any fixed angle $\vartheta$.
This amounts to solving 
\begin{equation}
\label{odetheta}
\frac{d}{d \sigma} z=e^{\mi \vartheta}f(z)
\end{equation}
on $Z=\mathbb{R}^{2N}$ by $\sigma\mapsto z(e^{\mi \vartheta}\sigma)$, locally for real time $\sigma$.
For example, we may decompose $t=\sigma+\mi \tau$ into real and imaginary parts.
Then the local flows $\Phi^\sigma$ and $\Phi^{\mi \tau}$ commute, for any argument $z_0$, as long as $\Phi^t(z_0)$ is defined on the rectangle spanned by $\sigma$ and $\mi \tau$ :
\begin{equation}
\label{flowc}
\Phi^\sigma\circ\Phi^{\mi\tau}= \Phi^{\mi\tau} \circ \Phi^\sigma.
\end{equation}
This follows from \eqref{odec} and Cauchy's theorem applied to $f\circ z$ on the rectangle.
Just as well, we may admit isolated blow-up points of $z(t)$ inside the rectangle, as long as $f\circ z$ does not generate residues there.
Examples are all meromorphic solutions $z(t)$ of \eqref{odec}, i.e. solutions which are complex analytic for all $t\in\mathbb{C}$, except at isolated poles of finite order.
More generally, isolated essential singularities are also covered, in absence of branching.
Indeed, the derivative $\dot z(t)$ is then residue-free anywhere.

Let us consider some planar examples.
The standard quadratic pendulum
\begin{equation}
\label{pend}
\ddot z-2 z^2+2=0,
\end{equation}
a special case of \eqref{pendulum},
can be written as a first order system in $(z,\dot z)$ with polynomial vector field $f$.
It possesses the meromorphic homoclinic orbit
\begin{equation}
\label{pendhom}
z(t)=\Gamma(t):=1-3/ \cosh^2(t)\,.
\end{equation}
The double poles of $\Gamma(t)$, and triple poles of $\dot\Gamma(t)$, are located at
\begin{equation}
\label{pendpoles}
t = \mathrm{i}\, (k+\tfrac{1}{2}) \pi\,,\  k\in\mathbb{Z}\,.
\end{equation}
Note the purely imaginary period $\pi\mi$ of the complexified homoclinic orbit.

With a scaling parameter $\lambda>0$, we also obtain homoclinic orbits $z(t):=\lambda^2\Gamma(\lambda t)$ of
\begin{equation}
\label{pendlambda}
\ddot z-2 z^2+2\lambda^4=0\,.
\end{equation}
The asymptotic equilibria $z=\pm\lambda^2$ and the homoclinic orbit $z(t)$ shrink to zero, for $\lambda\searrow 0$.
The poles at $t=\mi (k+\tfrac{1}{2})\pi/\lambda$, on the other hand, admit wider and wider horizontal $\eta$-strips \eqref{horstrip} of analyticity.
Under discretization, this provides stronger and stronger exponential splitting estimates \eqref{split}, but never ultra-exponential splitting.
With the same Hamiltonian core, this phenomenon will also scale homoclinic splittings \cite{Gelfreich03} near Bogdanov-Takens bifurcations \cite{ArnoldODE, GuckenheimerHolmes}, under discretizations.

Example \eqref{pend}, \eqref{pendhom} is a special case of the Weierstrass elliptic functions $\wp$ which satisfy the general  quadratic second order ODE
\begin{equation}
\label{wp}
\ddot\wp-6\wp^2+\tfrac{1}{2}g_2 =0\,.
\end{equation}
Indeed, conservation of energy $-\tfrac{1}{2}g_3$ reduces that second order ODE to the algebraic curve
\begin{equation}
\label{wpham}
\dot\wp^2 = 4\wp^3-g_2\wp-g_3\,;
\end{equation}
see for example \cite{Lang}.
For nonvanishing discriminant $\Delta:=g_2^3-27g_3^2$\,, the Weierstrass functions $\wp,\ \dot\wp$ generate the field of meromorphic, doubly periodic functions on a complex lattice associated to the modular Klein invariant $J=g_2^3/\Delta$.
The lattices, alias $J$, parametrize and identify the algebraic curves \eqref{wpham} as complex tori.
For real homoclinic orbits, i.e. for degenerate discriminant $\Delta=0$ and real coefficients $g_2^3=27g_3^2\neq 0$, the periodicity lattice degenerates to infinite ``real period'' and a nonzero imaginary period.

The above examples are scalar Hamiltonian pendulum equations.
\emph{Time reversible} second order ODEs, however, are not standard Hamiltonian in general.
Here, time reversibility just requires $z(-t)$ to be a solution whenever $z(t)$ is.
For some background on time reversibility techniques see for example \cite{Sevryuk,FiedlerHeinze,FiedlerTuraev,VanderbauwhedeFiedler}.
Reversible systems in one degree of freedom still possess first integrals, i.e.\ nontrivial conserved quantities $H=H(z,\tfrac{1}{2}\dot{z}^2)$.
In other words, $\dot H=0$ along solutions $z(t)$.

Scalar ODEs of the form $\ddot{z}+\dot{z}^2+\mathbf{q}(z)=0$ provide examples which are special cases of \eqref{pendulumm}.
To construct an explicit first integral, define $h(a,z)$ by integration of the inhomogeneous linear equation $\tfrac{\partial}{\partial z}h=-2h+\mathbf{q}$ with initial condition $h(a,a)=0$.
Then $H(z,\tfrac{1}{2}\dot{z}^2):=\tfrac{1}{2}\dot{z}^2+h(a,z)\equiv0$ is a first integral.
The real parameter $a$ identifies the intersection $(a,0)$ of the reversible orbit $(z,\dot z)$ with the horizontal $z$-axis.

To be specific, consider the time reversible, second-order, scalar ODE
\begin{equation}
\label{rev}
\ddot z +\dot z^2 +z^2-3z=0
\end{equation}
with quadratic nonlinearity $\mathbf{q}(z):=z^2-3z$.
By integrability $H\equiv 0$, the homoclinic orbit $z(t)=\Gamma(t)$ to the hyperbolic trivial equilibrium $z=0$, at parameter $a=0$, satisfies
\begin{equation}
\label{revint}
\tfrac{1}{2}\dot{z}^2=\exp(-2z)-1+2z\,-\,\tfrac{1}{2}z^2\,.
\end{equation}
In particular, $\Gamma(t)$ encounters non-meromorphic singularities in complex time $t$.
Indeed, theorem \ref{thmhom} implies that $\Gamma$ cannot be complex entire.
The term $\exp(-2z)$ in ODE \eqref{revint}, on the other hand, prevents meromorphic singularities.
Even if the unknown singularity is isolated and without branching, it will at least have to be essential.

In the scalar Hamiltonian pendulum examples \eqref{pend}, \eqref{wp} above, all nontrivial periodic orbits are also singular (meromorphic, by Weierstrass \cite{Lang}), rather than complex entire -- just like their homoclinic limit $\Gamma$ itself.
Our time reversible example \eqref{rev}, however, violates the Wittich condition \cite{Wittichm} which prevents non-polynomial entire solutions a priori.
In the language of Wittich, in fact, each of the \emph{two} terms $\dot z^2$ and $z^2$ is ``of maximal dimension'' $d=2$, where only a single such term is permitted.

Notably, \eqref{rev} does possess an explicit \emph{complex entire $2\pi$-periodic orbit} $\gamma$\,:
\begin{equation}
\label{revper}
z(t) = \gamma(t) := 2+\sqrt{2}\cos t\,.
\end{equation}
Under discretization, such complex entire periodic orbits will provide devil's staircases with ultra-exponentially thin resonance plateaus, and ultra-sharp Arnold tongues; see section \ref{UsA} and theorem \ref{thmultra-exp}.

After these illustrating examples, we now address \emph{local unstable manifolds} $W^u$ in our complex time setting \eqref{odec}.
The case of \emph{local stable manifolds} $W^s$ is analogous, in reverse time $t\mapsto -t$.

Let $W^u$ denote the set of all initial conditions $z_0$\,, such that $z(\sigma)$ remains in a suitably chosen small open neighborhood of the hyperbolic equilibrium $z=0$, for all real history $\sigma\leq0$.
By definition, the set $W^u$ is \emph{globally invariant} under the local flow $\Phi^\sigma$ on $Z=\mathbb{R}^{2N}$ in \emph{backward} real time $\sigma\leq 0$.

We now describe the (backward) invariant set $W^u$ as a local manifold.
Let $A:=f'(0)$ denote the linearization at the trivial equilibrium  $z=0$ of \eqref{odec}, i.e.
\begin{equation}
\label{odeclin}
\begin{aligned}
\dot z &= f(z) = Az+\ldots \,;\\
z(0) &= z_0\,.
\end{aligned}
\end{equation}
We assume $A$ is hyperbolic, i.e. all eigenvalues $\mu$ of $A$ possess nonzero real part $\Re\,\mu$.
Linear spectral decomposition $Z=Z^-\oplus Z^+$ into the generalized eigenspaces of eigenvalues with $\pm\,\Re\,\mu>0$, respectively, provides coordinates $z=(z^-,z^+)$ which split  \eqref{odeclin} as
\begin{equation}
\label{odecpm}
\begin{aligned}
\dot z^- &= A^-z^- + f^-(z^-,z^+)\,;\\
\dot z^+ &= A^+z^+ + f^+(z^-,z^+)\,.
\end{aligned}
\end{equation}
Note $\pm\,\Re\,\mathrm{spec}\,A^\pm >0$, and $f^\pm: Z\rightarrow Z^\pm$ are of at least quadratic order at $z=0$.

Standard theory then identifies $W^u$ as the graph of a local function
\begin{equation}
\label{wu}
U: Z^+_\mathrm{loc} \rightarrow Z^-_\mathrm{loc}\,,
\end{equation}
 in a small neighborhood $Z^-_\mathrm{loc} \oplus Z^+_\mathrm{loc} $ of $z=0$.
Therefore, $W^u$ is called the local unstable manifold of the hyperbolic equilibrium $z=0$.
See for example \cite{ChowHale}.
Here $U(0)=0$ marks the trivial equilibrium $z=0$.
The unstable eigenspace $Z^+$ is the tangent space to $W^u$ at $z=0$, i.e. $U'(0)=0$.

Perron's proof of the above result is based on the variation-of-constants formula, componentwise for \eqref{odecpm}, with given initial unstable component $z^+_0$ at $\sigma=0$.
The bounded initial condition for the stable component $z^-(\sigma_0)$ gracefully disappears exponentially from the resulting equations, for $\sigma_0\searrow -\infty$.
The past histories $z(\sigma)=z(\sigma;z^+_0)$, in the Banach space of bounded continuous functions for $\sigma\leq 0$, then result from a straightforward application of the implicit function theorem, with parameter $z^+_0$.
To obtain the above results, it only remains to define
\begin{equation}
\label{wudef}
U(z^+_0) := z^-(\sigma;z^+_0)\vert_{\sigma=0}\ .
\end{equation}
Except for the trick of an initial time $\sigma_0\searrow-\infty$ in the first equation of \eqref{odecpm}, the procedure is actually quite similar to existence and uniqueness proofs in semigroup theory; see for example \cite{Henry, Pazy}.

Based on the implicit function theorem, the procedure works just the same in complex notation for $Z=\mathbb{C}^N$, as long as we keep time $t=\sigma$ real.
Differentiability of $U$ however, as provided by the implicit function theorem, can now be viewed as \emph{complex differentiability}.
This proves that $U$ is locally holomorphic, in our case.
Local invariance of the local unstable manifold $W^u$ in \emph{imaginary time} $t=\mi\tau$ follows because the real and imaginary flows commute, by \eqref{flowc}, and because the implicit function theorem guarantees uniqueness.
If necessary, we may have to slightly reduce the open neighborhoods $Z^\pm_\mathrm{loc}$\,, for this argument.
As for invariance, in fact, imaginary time $\mi\tau$ may produce overflow at the boundary of $W^u$, just as positive real time $\sigma>0$ does.
(Poincaré linearization for nonresonant real spectra of $A^\pm$ will greatly facilitate this aspect, in the next section.)
In summary we obtain a reduced ODE
\begin{equation}
\label{odec+}
\begin{aligned}
\dot z^+&=f(U(z^+),z^+)=A^+z^++\ldots\ ,\\
z^+(0)&=z^+_0 \in Z^+_\mathrm{loc}\ ,
\end{aligned}
\end{equation}
for the complex analytic local flow $\Phi^t$ on the holomorphic local unstable manifold $W^u=\mathrm{graph}(U)\subseteq Z$, in complex time $t=\sigma+\mi\tau$.

\section{Nonresonant real spectra and imaginary quasi-periodicity} \label{Liq}

Dropping superscripts $^\pm$ temporarily, to simplify notation, we discuss local flows $\Phi^t$ of \eqref{odec+} in complex time $t=\sigma+\mi\tau$, rewritten as
\begin{equation}
\label{ode+}
\begin{aligned}
\dot z&=Az+g(z)\ ,\\
z(0)&=z_0 \in Z_\mathrm{loc}\ ,
\end{aligned}
\end{equation}
with $g(0)=0,\ g'(0)=0$.
In addition to local complex analyticity we assume unstable spectrum $\Re\,\mu_j>0$ for all eigenvalues $\mu_j\in\mathrm{spec}\,A$.

Assume $(Az)_j=\mu_jz_j$ is diagonal, and spectrally nonresonant in the sense of \eqref{nonres+}.
Then \emph{Poincaré linearization} at the spectrally nonresonant equilibrium $z=0$ implies that there exists a local, holomorphic, near-identity diffeomorphism
\begin{equation}
\label{pdif}
z=\phi(\zeta)=\zeta+\ldots
\end{equation}
which exactly linearizes \eqref{ode+}, locally, to become
\begin{equation}
\label{plin}
\dot \zeta_j=\mu_j\,\zeta_j\,,
\end{equation}
for $j=1,\ldots,N$. 
See for example \cite{ArnoldODE, Ushiki-N}.
For complete proofs and the definitive generalization including resonances, better see \cite{Ilyashenko}.
Differently from Ushiki, we cannot assume global analyticity of, say, the time-$\eps$ map of the underlying ODE flow, here.
In fact, the complex transformation $\phi(\zeta)=\zeta+\ldots$ need only be locally holomorphic.

Now assume that all eigenvalues $\mu_j$ are real.
Then the diagonal flow of \eqref{plin} in complex time $t=\sigma+\mi\tau$ decomposes as
\begin{equation}
\label{zetat}
\zeta_j(t) = e^{\mu_jt} \zeta_{j,0} = e^{\mu_j\sigma} \ e^{\mi\mu_j\tau}\ \zeta_{j,0}\,.
\end{equation}
In particular, the absolute values $|\zeta_j(t)|=e^{\mu_j\sigma}\,|\zeta_{j,0}|$ are independent of $\tau=\Im\,t$.
For $\sigma\rightarrow -\infty$, they decrease to zero, exponentially in $\sigma$ and uniformly in $\tau$.
Therefore
\begin{equation}
\label{qper}
z(\sigma+\mi\tau)=\phi\Big(e^{\mu_j\sigma} \ e^{\mi\mu_1\tau}\ \zeta_{1,0}\,,\,\ldots\,,\,e^{\mu_N\sigma} \ e^{\mi\mu_N\tau}\ \zeta_{N,0}\Big)
\end{equation}
is well-defined, for all sufficiently small $|z_0|$ and all $\tau\in\mathbb{R},\ \sigma\leq 0$.
In fact, the dynamics \eqref{zetat} of the argument $\zeta(t)$ of the analytic local diffeomorphism $\phi$ in \eqref{qper} is a parallel irrational flow on a torus of dimension at most $N$, in the imaginary time direction $\tau$.
Therefore the dynamics \eqref{qper} is quasi-periodic in $\tau$, for small $|z|$; see lemma \ref{lemap}.

\section{Almost-periodic Fourier series} \label{Fos}

For later use, we briefly recall some facts on almost-periodic functions, a slight generalization of the quasi-periodic functions which we just encountered. We follow the approach by Harald Bohr \cite{Bohr}, based on continuity and sup-norms. For $L^p$ variants and many others see also \cite{Besicovitch, Corduneanu, Haraux}.
As an important example we establish almost-periodicity of solutions in stable and unstable manifolds of equilibria $v$, in lemma \ref{lemap}, provided the associated tangential spectrum is real and nonresonant.  
 
Instead of Bohr, we start from the equivalent modern definition of almost-periodicity by Haraux \cite{Haraux}. 
Let $\mathcal{C}$ denote the Banach space of bounded continuous functions $z: \mathbb{R}\rightarrow Z=\mathbb{C}^N$, with the usual sup-norm.
Let $(\mathbf{T}_\theta z)(\tau):=z(\tau+\theta)$ denote the translates of $z$.
Then $z\in\mathcal{C}$ is called \emph{almost-periodic} if the set of all translates $\mathbf{T}_\theta z,\ \theta\in\mathbb{R},$ is relatively compact in $\mathcal{C}$.
See also \cite{Bohr}, §§ 42-44.
An example is the Riemann $\zeta$-function $z(\tau):=\sum_{n\geq1} n^{-(\sigma+\mi\tau)}$, which is almost-periodic in the imaginary direction $\tau$ with frequencies $\omega=-\log n$.

Let $z$ be almost-periodic.
Following \cite{Bohr}, § 56, we define the \emph{Fourier coefficient} $a_\omega$ of $z$ as the limiting average
\begin{equation}
\label{aom}
a_\omega \ := \ \lim_{T\rightarrow\infty} \ \frac{1}{T} \int_0^T e^{-\mi\omega\tau}z(\tau)\,d\tau\,.
\end{equation}
By almost-periodicity, these limits exist for all real frequencies $\omega$.
They are nonzero for at most countably many $\omega$.
Moreover, the Fourier coefficients $a_\omega$ are square summable: the finite limiting average of $|f|^2$ concides with $\sum_\omega |a_\omega|^2 < \infty$.

In this context, an almost-periodic function is called \emph{quasi-periodic} when the frequencies $\omega$ of nonzero Fourier coefficients $a_\omega$ are contained in a finitely generated $\mathbb{Z}$-module.
Periodic functions possess the frequency $\omega=2\pi/T$ of the minimal period $T$, as a single generator.


For any countable collection of nonzero $a_\omega\in\mathbb{C}, \ \omega\in\mathbb{R}$, consider formal Fourier series
\begin{equation}
\label{ser}
z(\tau) \ := \ \sum_\omega \,a_\omega \,e^{\mi\omega\tau}\,.
\end{equation}
The closure of all finite trigonometric sums \eqref{ser}, in the Banach space $\mathcal{C}$ of uniform convergence, are precisely the almost-periodic functions.
Conversely, any almost-periodic function $z(\tau)$ can be approximated, arbitrarily well and uniformly, by certain partial sums of its Fourier series \eqref{ser} with Fourier coefficients $a_\omega$ defined by \eqref{aom}.
Moreover, the representation \eqref{ser} of $z$ by any formal Fourier series is unique.
See \cite{Bohr}, §§ 84-92. 
For absolute convergence of the series see also § 93.

As an important example, we now address solutions $z(t)$ in the local unstable or stable manifold of a hyperbolic equilibrium $v$.
To be specific, we address the local unstable manifold $W^u$; for stable manifolds $W^s$ we just reverse time.
Notationally, we return to \eqref{ode+} in section \ref{Liq}.
In particular we assume the spectrum of the diagonalizable linearization $A$ at $v=0$ to be semisimple, unstable, and real nonresonant as in \eqref{nonres+}.

\begin{lem}\label{lemap}
In the real nonresonant spectral setting of \eqref{ode+}, \eqref{nonres+} above, consider the solution $z(t)$ of of any $z(0)\in W^u$, in complex time $t=\sigma+\mi\tau$.\\
Then $z(t)\in W^u$, for all $\sigma=\Re\,t\leq 0$.
Moreover $\tilde z(\tau):=z(\sigma+\mi\tau)$ is quasi-periodic in the imaginary time direction $\tau\in\mathbb{R}$.
The Fourier series of $\tilde z$ converges absolutely, and uniquely represents $\tilde z:$\begin{equation}
\label{gamser}
\begin{aligned}
\tilde z(\tau) &\phantom{:}= \ \sum_\omega \,a_\omega(\sigma) \,e^{\mi\omega\tau}\,,\ \mathrm{for}\\
a_\omega(\sigma) \,&:= \ \lim_{T\rightarrow\infty} \ \frac{1}{T} \int_0^T e^{-\mi\omega\tau}z(\sigma+\mi\tau)\,d\tau\,.
\end{aligned}
\end{equation}
In forward time $\sigma\geq 0$, and for nonresonance condition \eqref{nonres-}, the analogous statements hold for $z(0)\in W^s$.
\end{lem}

\begin{proof}
Without loss of generality, it is sufficient to study the case $z(\sigma)\rightarrow v$ for $\sigma\rightarrow -\infty$, with nonresonant unstable linearization $A^+$ at $v=0$.

First, consider $\tau=0$, so that $t=\sigma+\mi\tau=\sigma$ is real.
Then $z(\sigma) := \Gamma(\sigma)$ enters the local unstable manifold $W^u$ of $v=0$, and remains there, for all sufficiently large negative times $\sigma$; see \eqref{odeclin}--\eqref{odec+}. 

In particular, Poincaré linearization \eqref{ode+}--\eqref{plin}, applied to the reduced equation \eqref{odec+} in the local unstable manifold $W^u$, implies
\begin{equation}
\label{zsig}
z(t)=\big((U\circ\phi)(\zeta(t)), \phi(\zeta(t)\big),
\end{equation}
with $U: Z^+_\mathrm{loc}\rightarrow Z^-_\mathrm{loc}$ and a near-identity diffeomorphism $\phi$, both holomorphic near $0\in Z^+$.

Since the flow \eqref{plin}, \eqref{zetat} of $\zeta(t)$ is linear, and all unstable eigenvalues at $v=0$ are real, we can extend \eqref{zsig} to all imaginary times $\tau\in\mathbb{R}$, globally; see \eqref{qper} and the discussion of $W^u$-invariance in sections \ref{Cot} and \ref{Liq}.

The right hand side of \eqref{zsig} is holomorphic in $\zeta\in Z^+_\mathrm{loc}$\,.
We expand the local power series with respect to $\zeta$ and substitute $\zeta(t)$ in \eqref{zsig} by \eqref{zetat}, trigonometrically, as in \eqref{qper}.
This provides an explicit Fourier series \eqref{gamser} for $\tilde{z}(\tau):=z(\sigma+\mi\tau)$.

By geometric series majorants, inside the open disc of local analyticity, the Fourier series \eqref{gamser} converges absolutely and represents $\Gamma$.
In particular, $\Gamma(\sigma+\mi\tau)$ is almost-periodic in imaginary time $\tau$, as claimed.
By local analyticity of $\phi$ and $U$, the frequency module is finitely generated by $\mathrm{spec}\, A^+$.
This proves the lemma.
\end{proof}

\section{Fourier invariants} \label{Foi}

We can now apply the results of the previous sections to homo- or heteroclinic orbits $\Gamma(t),\ t=\sigma+\mi\tau$, between asymptotic hyperbolic equilibria 
\begin{equation}
\label{vpm}
v^\pm \ =\ \lim_{\sigma\rightarrow\pm\infty} \Gamma(\sigma+\mi\tau);
\end{equation}
see \eqref{homhet}. 
At $v^\pm$ we assumed diagonalizable $A^\mp$ with semisimple, real, nonresonant spectra; see \eqref{nonres-}, \eqref{nonres+}.

Lemma \ref{lemap} has established almost- and quasi-periodicity of $\tilde z(\tau)=z(\sigma+\mi\tau):=\Gamma(\sigma+\mi\tau)$, in imaginary time $\tau$, for all fixed and sufficiently large $|\sigma|$. 
Indeed, $\tilde z(0)\in W^u$ and $\tilde z(0)\in W^s$ then hold, respectively and separately for each sign of $\sigma$. 
Let $a_\omega(\sigma)$ denote the Fourier coefficients of the almost-periodic functions $\tilde z(\tau)$, as in \eqref{gamser}.
For complex entire $\Gamma(t)$, we now establish that suitably scaled Fourier coefficients $c_\omega$
neither depend on $\sigma$, nor its sign, for all  sufficiently large $|\sigma|$.
By contradiction, this fact will prove theorems \ref{thmhom} and \ref{thmhet} in section \ref{Pro}.

\begin{lem}\label{leminv}
Suppose $\Gamma(t)$ is a complex entire homo- or heteroclinic orbit \eqref{vpm} between equilibria $v^\pm$ with separately nonresonant real spectra, as in \eqref{nonres-}, \eqref{nonres+} above.
Fix any real frequency $\omega$. \\
Then the Fourier coefficients $a_\omega(\sigma)$ of the almost-periodic functions $\tilde z(\tau):=\Gamma(\sigma+\mi\tau)$ produce $\sigma$-invariant scaled constants
\begin{equation}
\label{cominv}
c_\omega := e^{-\omega\sigma}a_\omega(\sigma),
\end{equation}
for all sufficiently large $|\sigma|$, independently of $\sigma$ and independently of the sign of $\sigma$.
\end{lem}

\begin{figure}[t]
\centering \includegraphics[width=0.9\textwidth]{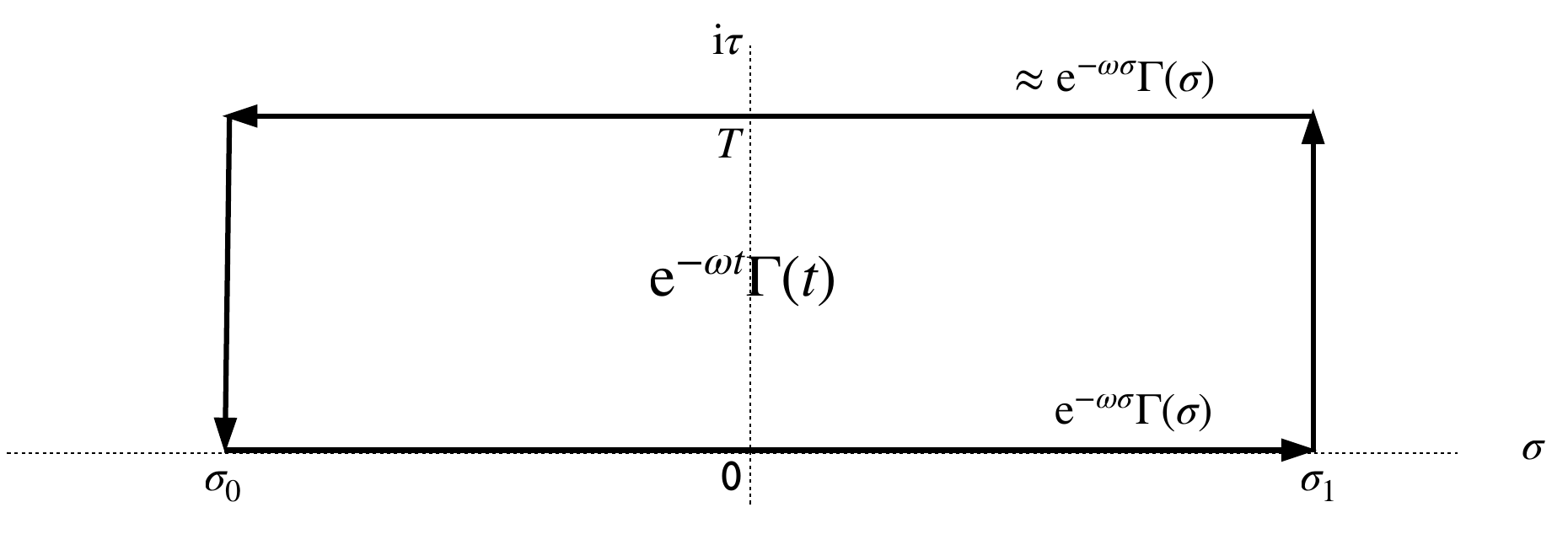}
\caption{\emph{
Cauchy theorem along a left-oriented rectangle in complex time $t=\sigma+\mi\tau$, with corners $\{\sigma_0,\sigma_1,\sigma_1+\mi T,\sigma_0+\mi T\}$.
By almost-periodicity, $T$ can be chosen such that the integrands $\tfrac{1}{T} e^{-\omega t}\Gamma(t)$ on the top and bottom horizontal boundaries remain uniformly $\delta $-close; see \eqref{estint}.
Therefore the limits of the vertical averages $c_\omega$  of lemma \ref{leminv} coincide, in \eqref{cominv}, at the left and right boundaries $\Re\,t = \sigma_0\,,\sigma_1$\,.
}}
\label{fig2}
\end{figure}

\begin{proof}
By definition \eqref{aom} of the Fourier coefficients $a_\omega$\,, we have to determine
\begin{equation}
\label{comint}
\begin{aligned}
c_\omega := e^{-\omega\sigma}\, a_\omega(\sigma) &= \lim_{T\rightarrow\infty} \frac{1}{T} \int_{0}^{T} e^{-\omega\sigma}e^{-\mi\omega\tau}\,\Gamma(\sigma+\mi\tau)\,d\tau =\\
&=  \lim_{T\rightarrow\infty} \frac{-\mi}{T} \int_{\sigma}^{\sigma+\mi T} e^{-\omega t} \,\Gamma(t)\, dt\,,
\end{aligned}
\end{equation}
with $t:=\sigma+\mi\tau$, as always.
To show independence of $c_\omega$ on $\sigma=\Re\,t$, we will invoke Cauchy's theorem for the analytic integrand of the second line. See fig. \ref{fig2}.

First, fix real $\sigma_0,\ \sigma_1$ of large absolute value, but otherwise arbitrarily and, in particular, of arbitrary sign.
Next, consider $T>0$ large enough such that the limits in \eqref{comint} are achieved, at $\sigma=\sigma_0,\sigma_1$, up to some arbitrarily small error $\eps>0$.
By almost-periodicity of the pair $(e^{-\mi \omega\tau},\ \Gamma(\sigma_0+\mi \tau))$ in $\tau=\Im\,t$, the values of $e^{-\mi \omega T}$ and $\Gamma(\sigma_0+\mi T)$ simultaneously get arbitrarily close to $1$ and $\Gamma(\sigma_0)$, again and again, for $T\rightarrow +\infty$.
By continuous dependence of the complex flow \eqref{odec} on initial data, we may therefore choose $T>0$ large enough, such that in addition
\begin{equation}
\label{estint}
 \frac{1}{T} e^{-\omega\sigma} \,\vert e^{-\omega\mi T}\Gamma(\sigma+\mi T)-\Gamma(\sigma)\vert\leq \delta\,,
\end{equation}
uniformly for $\sigma$ in the compact interval between fixed $\sigma_0$ and $\sigma_1$\,.

We can now apply Cauchy's theorem to the complex entire integrand $e^{-\omega t} \,\Gamma(t)$ on the rectangle with corners $t\in\{\sigma_0,\sigma_1,\sigma_1+\mi T,\sigma_0+\mi T\}$.
By construction, the desired difference of the integrals along the vertical boundaries equals the difference along the horizontal boundaries.
Therefore estimate \eqref{estint} of the integrands implies 
\begin{equation}
\label{estinv}
|c_\omega(\sigma_1)-c_\omega(\sigma_0)| \leq 2 \eps + |\sigma_1-\sigma_0| \delta < 3\eps,
\end{equation}
for sufficiently small $\delta=\delta(\sigma_0,\sigma_1,\eps)$.
Since $\eps>0$ was arbitrary, this proves the lemma.
\end{proof}

One subtlety concerns the Fourier coefficients $a_\omega(\sigma)$ in \eqref{gamser}. 
An easy adaptation of the arguments in the proof of lemma \ref{leminv} shows that the limits $a_\omega(\sigma)$ in fact exist, for any fixed real $\sigma$, and invariance \eqref{cominv} remains valid, as long as the Cauchy theorem applies.
Convergence of the Fourier series \eqref{gamser}, however, may fail, unless $|\sigma|$ is large enough.
The reversible pendulum \eqref{pend}, \eqref{pendhom} is an instructive exercise, where $c_\omega(\sigma) = c_{-\omega}(-\sigma)$ for nonzero $\sigma=\Re\, t$.

\section{Proof of main results} \label{Pro}

In this section we prove theorems \ref{thmhom} and \ref{thmhet}, by contradiction.
Therefore, contrary to the claims, we assume the existence of a complex entire homo- or heteroclinic orbit $\Gamma(t),\ t=\sigma+\mi\tau$, between hyperbolic equilibria $v^\pm$, as in the previous section.
In particular, the spectra of $A^\mp$ are assumed to be semisimple, real and nonresonant as in \eqref{nonres-}, \eqref{nonres+}.

\emph{Proof of theorem \ref{thmhet}.}
Our starting point is the Fourier series \eqref{gamser} of $\Gamma$ in imaginary time $\tau$.
For complex entire heteroclinic orbits $\Gamma(t)$ between hyperbolic equilibria $v^-\neq v^+$, we simply evaluate the averages \eqref{cominv}, \eqref{gamser} for $\omega=0$\,:
\begin{equation}
\label{c0vpm}
c_0 = \lim_{\sigma\rightarrow\pm\infty} a_0(\sigma) = \lim_{\sigma\rightarrow\pm\infty} \lim_{T\rightarrow\infty} \ \frac{1}{T} \int_0^T \Gamma(\sigma+\mi\tau)\,d\tau = v^\pm\,.
\end{equation}
To interchange the limits, for $\sigma\rightarrow -\infty$, we have used $\tau$-uniform convergence of $\Gamma(\sigma+\mi\tau) \in W^u(v^-)$ to the asymptotic equilibrium $v^-$.
This is due to nonresonance \eqref{nonres+} of the real unstable spectrum at $v^-$, and to Poincaré linearization \eqref{ode+}, \eqref{qper}.
For $\sigma\rightarrow +\infty$ and $v^+$, we argue analogously, in forward time.
Since \eqref{c0vpm} contradicts the heteroclinicity assumption $v^-\neq v^+$, the theorem is proved. 
\hfill $\bowtie$

\emph{Proof of theorem \ref{thmhom}.}
For homoclinic orbits $\Gamma$, we may assume $\Gamma\neq v^\pm=0$.
Lemmata \ref{lemap} and \ref{leminv} combine to show
\begin{equation}
\label{}
\Gamma(t) =  \ \sum_\omega \,a_\omega(\sigma)\,e^{\mi\omega\tau} \ = \ \sum_\omega \,c_\omega\,e^{\omega t}\,,
\end{equation}
for $t=\sigma+\mi\tau$ and all large enough $|\sigma|=|\Re\, t|$.
To obtain a contradiction, we will show
\begin{equation}
\label{com=0}
c_\omega = e^{-\omega\sigma}\, a_\omega(\sigma) = 0\,,
\end{equation}
for all real frequencies $\omega$.
This contradicts $\Gamma \neq 0$.

Consider $\omega=0$ first. 
Then \eqref{c0vpm} with $v^\pm=0$ implies \eqref{com=0}.

Consider frequencies $\omega<0$ next.
From the proof of theorem \ref{thmhet}, we recall $\tau$-uniform convergence of $\Gamma(\sigma+\mi\tau) \in W^u$ to $v^-=0$, for $\sigma\rightarrow -\infty$.
By nonresonance \eqref{nonres+} of the semisimple real unstable spectrum at $v=0$, and by \eqref{gamser}, we obtain uniform boundedness of the Fourier coefficients $a_\omega(\sigma)$, for $\sigma\rightarrow -\infty$.
In \eqref{cominv}, that same limit proves \eqref{com=0} for $\omega<0$.

To prove \eqref{com=0} for $\omega>0$, finally, we just consider $\sigma\rightarrow +\infty$, and argue via the stable manifold and nonresonance \eqref{nonres-} of the real stable spectrum at the equilibrium $v^+=0$.

This proves claim \eqref{com=0}, for all frequencies $\omega$, and completes the proof of the theorem.
 \hfill $\bowtie$

\section{Ultra-sharp Arnold tongues}\label{UsA}

Theorems \ref{thmhom} and \ref{thmhet} can be read as hints where \emph{not} to look for ultra-exponential separatrix splittings and ultra-invisible chaos.
Corollary \ref{corblowup} reveals how real-time homoclinic or heteroclinic orbits $\Gamma(t)$ are linked to blow-up in imaginary time, at least in the case of semisimple, real, separately nonresonant spectrum \eqref{nonres-}, \eqref{nonres+} at the asymptotic equilibria $v^\pm$.
Our explicit time reversible example \eqref{rev}, in contrast, features a complex entire periodic orbit $\gamma$\,; see \eqref{revper}.
We therefore comment on ultra-exponentially narrow plateaus of rational rotation numbers, alias ultra-invisible phaselocking or ultra-sharp Arnold tongues, in discretizations of such periodic orbits; see theorem \ref{thmultra-exp}.
For general background on discretizations of periodic orbits, see for example \cite{Beyn, Chenciner, Stoffer-multi, Stoffervar, Kloeden}, among others.

Our approach closely follows the treatment of exponentially small homoclinic splittings in \cite{FiedlerScheurle}, which we presented in detail long ago.
Along the way we repair a minor error in the distributed Poincaré section chosen there; see \eqref{Lker} below.

\subsection{Lyapunov-Schmidt reduction}\label{LS}
Adapting notation slightly, for this last section, let $\gamma(t),\ t\in\mathbb{C},$ be a nonstationary, complex entire orbit of minimal real period 1, for a complex entire ODE \eqref{odef}, at some real parameter $\lambda=0$.
As in \eqref{fcc}, we assume $f$ to be real, for real arguments.
In particular, the periodic orbit $\gamma(\sigma)$ is assumed real, for real times $\sigma$. 
As in \eqref{odeg}, we subsume $p$-th order one-step discretizations of real step size $\eps>0$ under time-$\eps$ stroboscope maps of non-autonomous ODEs
\begin{equation}
\label{odegc}
\dot{z}(\sigma)\ =\ f(\lambda,z(\sigma)) + \eps^p g(\lambda,\eps,\sigma/\eps,z(\sigma))\,,
\end{equation}
in real time $\sigma=\Re\,t$.
We assume $g$ to be complex entire in $\lambda,\eps,z$ and, analogously to  $f$ in \eqref{fcc}, real for real arguments.
We also assume explicit forcing period 1 of $g$, in the real variable $\beta=\sigma/\eps$.
The real parameter $\lambda$ will allow us to track solutions $z(\sigma)$ of fixed period 1, for $\eps>0$.
However, our technique will require commensurately resonant discretization steps
\begin{equation}
\label{epsn}
\eps=\eps_n:=1/n\,,
\end{equation}
with large enough $n\in\mathbb{N}$.

\begin{figure}[t]
\centering \includegraphics[width=0.6\textwidth]{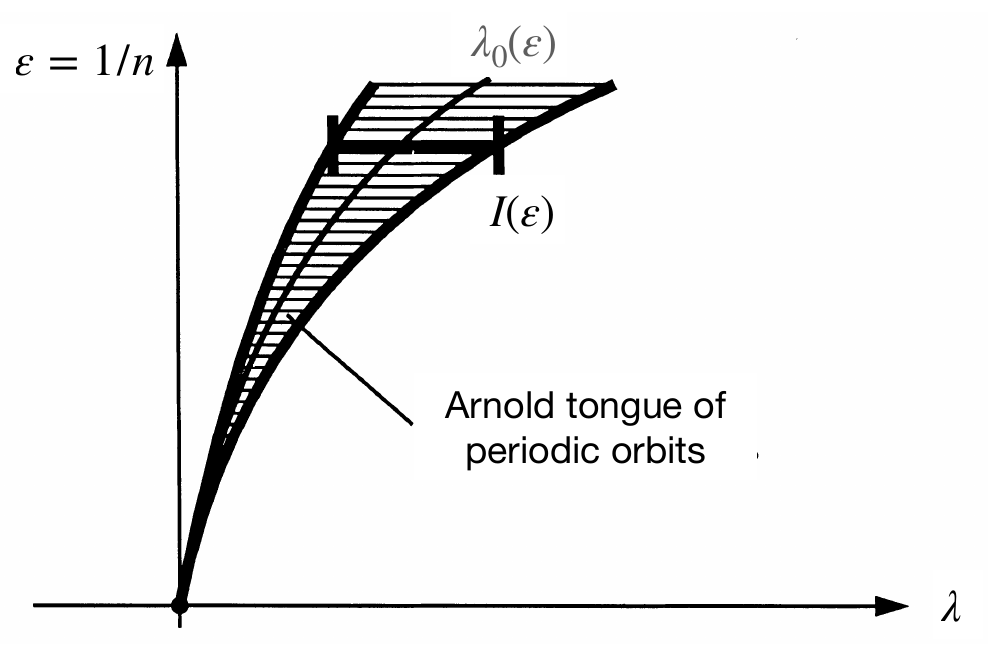}
\caption{\emph{
Schematic Arnold tongue for a 1-periodic orbit $\gamma$ of the ODE-flow \eqref{odef}, under discretization with step size $\eps$ or, equivalently, under rapid $\eps$-periodic non-autonomous forcing \eqref{odeg}.
Compare the homoclinic splitting region of fig. \ref{fig1}.
For the definition of $\lambda_0(\eps)$ see \eqref{avg}.
To work in spaces of 1-periodic functions, we now restrict to rapid non-autonomous forcings with commensurate discrete periods $\eps=1/n$.
For complex entire periodic orbits, theorem \ref{thmultra-exp} provides \emph{ultra-exponential upper estimates} for the width $\ell(\eps)$ of the real resonance (or phaselocking) intervals $I(\eps)$; see \eqref{ultraplat}.
In other words, the resulting Arnold tongue is \emph{ultra-sharp}.
For better visibility, the width of the Arnold tongue has been exaggerated.
}}
\label{fig3}
\end{figure}

To detect 1-periodic solutions, we now follow \cite{FiedlerScheurle} almost verbatim.
To keep track of complex time $t$ in $z$, and of real period 1 in $g$, we introduce two artificial time-shift parameters: $\alpha\in\mathbb{C}$ and $\beta\in\mathbb{R}/\mathbb{Z}$.
For real $\alpha,\beta$ this would actually be redundant; see \eqref{zred}, \eqref{Bred} below.
Complex $\alpha$, however, will keep track of complex analyticity in the reference periodic orbit $\gamma$.
Real $\beta$, in contrast, carry real periodicity of the rapid forcing $g$, and hence of discretization.
We then rewrite \eqref{odegc} for the deviation $\zeta(\sigma):=z(\sigma)-\gamma(\sigma+\alpha)$, in real time $\sigma$, but with complex parameter $\alpha$, as
\begin{eqnarray}
\label{L-F=0}
\mathcal{L}(\alpha)\zeta - \mathcal{F}(\lambda,\eps,\alpha,\beta,\zeta) & = & 0\,; \\
\label{Lker}
\int_{0}^{1} \zeta(\sigma)^T\cdot\dot{\gamma\,}(\sigma+\alpha)\,d\sigma & = & 0\,.
\end{eqnarray}

With the complex Floquet linearization $A(t):= \partial_z f(0,\gamma(t))$, we use the abbreviations
\begin{eqnarray}
\label{defL}
(\mathcal{L}(\alpha)\zeta)(\sigma) &:=& \dot\zeta(\sigma)-A(\sigma+\alpha)\zeta(\sigma)\,; \\
\label{defF}
\quad \mathcal{F}(\lambda,\eps,\alpha,\beta,\zeta)(\sigma) &:=& f(\lambda,\gamma(\sigma+\alpha)+\zeta(\sigma)) -  f(0,\gamma(\sigma+\alpha)) - A(\sigma+\alpha)\zeta(\sigma)\\
& & + \ \eps^p g(\lambda,\eps,\beta+\sigma/\eps,\gamma(\sigma+\alpha)+\zeta(\sigma))\,.\nonumber
\end{eqnarray}

We can then solve system \eqref{L-F=0}, \eqref{Lker} as a holomorphic map from the Banach space $C^1(\mathbb{R}/\mathbb{Z},\mathbb{C}^N)$ of complex-valued, real 1-periodic $\zeta$, to $C^0(\mathbb{R}/\mathbb{Z},\mathbb{C}^N)\times \mathbb{C}$.
A minor technical point concerns holomorphy with respect to, both, $\zeta$ and the complex time shift $\alpha$.
In the distributed version \eqref{Lker} of a Poincaré section, we use the standard scalar product $a^T\cdot b=a_1b_1+\ldots a_Nb_N$ for complex arguments. 
This replaces the Hermitian scalar product $\bar{a}^T\cdot b$, erroneously used in (3.4b) of \cite{FiedlerScheurle}.
Note the trivial equilibrium solution $\zeta=0$, for $\lambda=\eps=0$ and all $\alpha,\beta$.
Moreover, $\mathcal{L}(\alpha)$ is Fredholm of index zero.

Henceforth, we \emph{assume} the original periodic orbit $\gamma(t)$ to be hyperbolic, i.e. with algebraically simple trivial Floquet multiplier 1.
In particular, $\dot{\gamma\,}$ spans the trivial kernel of the Fredholm operator $\mathcal{L}(0)$.
As a second nondegeneracy assumption, we require that $\partial_\lambda f(0,\gamma)$ spans the one-dimensional co-kernel of $\mathcal{L}(0)$.
Analogously, by complex time-shift $\alpha$, we note that $\dot{\gamma\,}(\cdot+\alpha)$ and $\partial_\lambda f(0,\gamma(\cdot+\alpha))$ then span the kernel and co-kernel of $\mathcal{L}(\alpha)$, respectively.

We can now apply standard Lyapunov-Schmidt reduction, in the general spirit of \cite{Vanderbauwhede} or \cite{Wulff}. 
More specifically, see \cite{ChowHale}, chapter 11.
We skip quite a few technical details which are analogous to \cite{FiedlerScheurle}.
The modified Poincaré section \eqref{Lker} still eliminates the trivial kernel $\dot{\gamma\,}(\cdot+\alpha)$. 
Indeed, the integral \eqref{Lker}, evaluated at the kernel $\zeta(\sigma)=\dot{\gamma\,}(\sigma+\alpha)$, is independent of $\alpha$, e.g. by Cauchy's theorem.
At $\alpha=0$, i.e. along the real periodic orbit $\gamma(\sigma)$, the integral is nonzero positive.

Therefore, we obtain solutions $\zeta(\lambda,\eps,\alpha,\beta)(\sigma)$ from the trivial solution $\zeta(0,0,\alpha,\beta)=0$, for all time-shift parameters $\alpha, \beta$.
The redundancy of the time shifts $\alpha$ and $\beta$, in the real case, implies that
\begin{equation}
\label{zred}
\zeta(\lambda,\eps,\alpha,\beta)(\sigma+\sigma') = \zeta(\lambda,\eps,\alpha+\sigma',\beta+\sigma'/\eps)(\sigma) 
\end{equation}
holds for all real $\sigma, \sigma'$. 
The remaining complex one-dimensional co-kernel of $\mathrm{range}\, \mathcal{L}(\alpha)$ provides a complex scalar Lyapunov-Schmidt reduced bifurcation equation $B(\lambda,\eps,\alpha,\beta)=0$.
Here \eqref{zred} with $\sigma':=-\beta\eps$ induces the intriguing redundancy
\begin{equation}
\label{Bred}
B(\lambda,\eps,\alpha-\beta\eps,0)=B(\lambda,\eps,\alpha,\beta)=0\,.
\end{equation}
In other words, we may drop the last argument $\beta\in\mathbb{R}/\mathbb{Z}$ of $B$, altogether, and note that we only have to solve
\begin{equation}
\label{B=0}
B(\lambda,\eps,\alpha)=0\,.
\end{equation}
But now $B$, unlike $\zeta$ before, is analytic \textbf{and} of real period $\eps$, in $\alpha$.

\subsection{Exponentially thin resonance plateaus}\label{Plat}
Analyticity in \emph{fixed} horizontal strips $|\Im\ \alpha | \leq \eta$, jointly with real period $\eps$, are the source of exponential smallness and ``invisibility'' of resonance plateaus; see for example lemma 3.4 in \cite{FiedlerScheurle}, or appendix A in \cite{Gelfreich02}.
Analyticity in \emph{arbitrary} horizontal strips $|\Im\ \alpha | \leq \eta$, analogously, becomes the source of their ultra-exponential smallness and ``ultra-invisibility'' in theorem \ref{thmultra-exp} below.

Note the trivial solution $B(0,0,\alpha)=0$, for all complex time-shift parameters $\alpha$.
Since we assumed that $\partial_\lambda f(0,\gamma(\cdot+\alpha))$ spans the co-kernel of $\mathcal{L}(\alpha)$, we also obtain $\partial_\lambda B(\lambda,\eps,\alpha)\neq0$.
The implicit function theorem therefore produces locally unique 1-periodic orbits of system \eqref{L-F=0}, \eqref{Lker} at
\begin{equation}
\label{lam}
\lambda=\lambda(\eps,\alpha)\,.
\end{equation}
Here $\lambda$ is real, for real arguments $\eps,\alpha$.
However, we had to restrict to discrete values of small $\eps=\eps_n=1/n>0$; 
the case $\eps=0$ formally corresponds to $n=\infty$ in \eqref{epsn}. 
With respect to those small $\eps_n\geq 0$, we can therefore only assert continuity of $\lambda$ at $\eps=0$.
Still, we obtain the following result on the width $\ell(\eps)$ of the 1-periodic \emph{resonance plateau}, $I(\eps)=\lambda(\eps,\mathbb{R})$.

\begin{thm}\label{thmultra-exp}
Consider a real 1-periodic orbit $\gamma$, in the above setting and under the above two nondegeneracy assumptions.
Let $\eps=\eps_n=1/n$ denote sufficiently fine discretization, alias rapid periodic forcing \eqref{odegc}.\\
(i) Assume $\gamma(t)$ is analytic in a neighborhood of the closed complex horizontal strip $|\Im\ t|\leq\eta$, for some $\eta>0$.
Then the resonance plateaus $I(\eps)$ are exponentially thin.
In other words, there exists a constant $C_\eta$ such that for every sufficiently small $\eps=\eps_n=1/n>0$ the width $\ell(\eps)$ of the resonance plateau $I(\eps)$  satisfies an \emph{exponential upper estimate}
\begin{equation}
\label{ultraplat}
\ell(\eps)\leq C_\eta \exp(-\eta/\eps)\,.
\end{equation}
(ii) If $\gamma(t)$ is complex entire, then the resonance plateaus $I(\eps)$ are ultra-exponentially thin.
In other words, for every $\eta>0$ there exists a constant $C_\eta$ such that an \emph{ultra-exponential upper estimate} \eqref{ultraplat} holds, for every sufficiently small $\eps=\eps_n=1/n>0$.
\end{thm}

For a schematic illustration of an ultra-exponentially thin wedge \eqref{ultraplat} see fig. \ref{fig3}.
Differently from the homoclinic case of fig. \ref{fig1}, however, admissible values for $\eps$ are now restricted to discrete levels $\eps=1/n$\,.

\begin{proof}

With respect to $\eps=\eps_n\geq 0$\,, the implicit function theorem only provides continuity of $\lambda(\eps,\alpha)$, at $\eps=0$.
With respect to $\alpha$, however, we can invoke complex analyticity and real $\eps$-periodicity of $B$.
To prove claim (i), we now fix $\eta>0$ and consider $\alpha$ in the fundamental rectangle $0\leq \Re\,\alpha\leq\eps,\ |\Im\,\alpha|\leq\eta$.
In the complex entire case of claim (ii), we proceed with arbitrary $\eta>0$, instead.

By the procedure of Lyapunov-Schmidt reduction, $B$ is bounded on that rectangle, uniformly for sufficiently small $|\lambda|, \eps$. 
The same holds true for the complex partial derivative $\partial_\alpha B$.

Fourier-expansion with respect to $\Re\,\alpha$ then implies exponential smallness for $\partial_\alpha B$, as in \eqref{ultraplat},  uniformly for real $\alpha$ and small $|\lambda|, \eps$. 
The partial derivative $\partial_\lambda B\neq 0$, in contrast, remains bounded away from zero, uniformly in the same region.

Implicit differentiation then implies uniform exponential smallness of $\partial_\alpha\lambda=-\partial_\alpha B/\partial_\lambda B$, for real $\alpha$.
By real $\eps$-periodicity of $\lambda=\lambda(\eps,\alpha)$ in $\alpha$, we obtain the exponential estimate \eqref{ultraplat} for the length $\ell(\eps)$ of the resonance plateau $I(\eps)$, and the theorem is proved.
\end{proof}

Offhand, our approach only yields continuity, in $\eps$, of the ``average'' location
\begin{equation}
\label{avg}
\lambda_0(\eps) := \frac{1}{\eps} \int_0^\eps \lambda(\eps,\alpha)\,d\alpha
\end{equation}
in the resonance plateau, i.e. of the zero-th Fourier coefficient of the $\eps$-periodic parameter locations $\alpha\mapsto\lambda(\eps,\alpha)$.
See fig. \ref{fig3}.
Subject to the perturbation $\eps^pg$ , of course, this $\eps$-dependence is actually of order $\eps^p$\,; see also \cite{Beyn}.
The resulting resonance plateau $I(\eps)$, however, will be exponentially or ultra-exponentially thin, respectively.

\subsection{Sharp Arnold tongues}
For the ODE flow $\eps=0$ of \eqref{odegc} in real time $\sigma$, our two nondegeneracy assumptions provide a local path of hyperbolic periodic orbits $\sigma\mapsto\gamma=\gamma(\lambda,\sigma)$ with minimal period $T=T(\lambda)>0$ in $\sigma
$.
Under periodic forcings of small amplitude $\eps$, but say of fixed (non-rapid) period $T'>0$, the path turns into a local path of invariant closed curves which remain normally hyperbolic.
Let $\rho(\lambda,\eps)$ denote the rotation numbers of the resulting time-$T'$ circle diffeomorphism \cite{Hartman}.
For example, $\rho(\lambda,0)=T'/T(\lambda)$ in the ODE flow case.
For generic one-parameter families of smooth diffeomorphisms, e.g. for constant $\eps$, the rotation numbers $\rho(\lambda,\eps)$ turn out to be locally constant, independently of $\lambda$, at any rational value of $\rho$. 
This phenomenon is called \emph{devil's staircase}.
See for example \cite{Brunovsky, Devaney, Ligeti}, and others.

Let us therefore fix a rational value $\rho=m_0/n_0$.
The \emph{Arnold tongue} phenomenon then describes how the resonance zones, i.e. the preimages $\rho^{-1}(m_0/n_0)$ in the $(\lambda,\eps)$-plane, emanate from their tips at $(\lambda,0)$, with $\lambda$ given by $T(\lambda)=T'n_0/m_0$.
For all small forcing amplitudes $\eps>0$, but still non-rapid fixed forcing period $T'>0$, the Arnold tongue provides resonance plateaus of width $\ell(\eps)$.
Generically, in fact these tips turn out to be only algebraically sharp, rather than (ultra-)exponential: typically, $\ell(\eps)$ is of order $\eps^{n_0/2-1}$, for $n_0\geq5$.

In our approach, we have technically been limited to diffeomorphisms of step size $\eps=1/n$.
Consider therefore sequences $n=a n_0\,,\ m=am_0$\,, for some integer scaling parameter $a \nearrow \infty$.
Then the $m$-th iterate of our stroboscopic $\eps$-discretizations \eqref{odegc}, now with rapid forcing period $\eps=T(0)/n$, corresponds to the above case with fixed $T'=m\eps=T(0) m_0/n_0$\,, i.e. to the fixed rotation number $\rho=m_0/n_0$\,.
Although we are limited to small discrete values $\eps=\eps_{an_0}$ at integer $a\nearrow\infty$, theorem \ref{thmultra-exp} constrains these ultra-invisible horizontal resonance plateau sections of the Arnold tongue  to an \hbox{(ultra-)exponentially} sharp wedge region.
This does not contradict the previous generic result, because our nongeneric time-$T'$ maps now possess the time-$\eps$ stroboscope map, as an $m$-th root; see \eqref{root} in section \ref{Root}.
In conclusion, Arnold tongues caused by discretization of complex entire periodic orbits, like example \eqref{rev}, \eqref{revper}, turn out ultra-sharp.

\subsection{Hamiltonian and reversible pendula}
Theorem \ref{thmultra-exp} does not cover Hamiltonian and reversible systems, directly.
Indeed, the trivial Floquet multiplier 1 is at least algebraically double in those cases, because families of periodic orbits typically occur in paths, even in absence of further parameters $\lambda$.
One remedy is an adaptation of our Lyapunov-Schmidt reduction process to that situation; see for example the techniques and references in \cite{FiedlerHeinze,VanderbauwhedeFiedler}.
As a quick alternative, we briefly sketch a one-parameter embedding of periodic and homoclinic orbits for reversible or Hamiltonian pendulum equations of the scalar form
\begin{equation}
\label{grev}
\ddot z + \mathbf{q}(z,p)=0,\qquad \textrm{where\ }p:=\tfrac{1}{2}\dot z^2;
\end{equation}
see also example \eqref{rev}. \eqref{revint}.
For simplicity, we start from a real first integral of the general form
\begin{equation}
\label{Hrev}
H=H(z,p)= \textrm{const}.
\end{equation}
with partial derivative $H_p>0$.
In other words, we study nonlinearities $\mathbf{q}=H_z/H_p$\,.
Let us then consider damped second order ODEs of the specific form
\begin{equation}
\label{Hlambda}
\ddot z + \frac{1}{H_p}(H-\lambda)\dot z + \mathbf{q} =0,
\end{equation}
with real parameter $\lambda$.
Then $H$ relaxes to the level curve $\{H=\lambda\}$, in the real $(z,\dot z)$-plane, except at equilibria $\mathbf{q}(z,0)=0$. 
Indeed,
\begin{equation}
\label{Hdot}
\dot H=(\lambda-H)\dot z^2
\end{equation}
along any solution of \eqref{Hlambda}.
In particular, this embeds real periodic orbits $\gamma$ of \eqref{grev} with first integral \eqref{Hrev} into the parametrized ODE \eqref{Hlambda}, for real parameters $\lambda$ given by the constant value of $H$ on $\gamma$.
In the analytic case, complex time extensions follow suit.
Theorem \ref{thmultra-exp} and the above remarks then apply to the setting of \eqref{Hlambda}.

For the second order pendulum, the nonlinearity $\mathbf{q}=\mathbf{q}(z)$ is independent of $p=\tfrac{1}{2}\dot z^2$.
The Hamiltonian integral $H$ is then given by $H(z,p)=p+\mathbf{Q}(z)$, where the potential $\mathbf{Q}$ is a primitive of $\mathbf{q}$.
In that case, \eqref{Hlambda} simplifies to the damped version
\begin{equation}
\label{Hdamped}
\ddot z + (H-\lambda)\dot z + \mathbf{q} =0.
\end{equation}

Alas, our results have only provided an ultra-exponential upper estimate for Arnold tongues associated to complex entire periodic orbits. 
This raises the question: how ultra-sharp are they, really?
Complementary lower estimates of ultra-sharp Arnold tongues are not known.
For homoclinic orbits, we are even lacking \emph{any} 1,000\ \euro\ complex entire example with ultra-exponential upper splitting estimate.


\end{document}